\documentclass[a4paper,11pt]{article}
\pdfoutput=1

\usepackage{jheppub}
\usepackage{amssymb}
\usepackage{amsmath}
\usepackage{mathdots}
\usepackage{verbatim}
\usepackage[utf8]{inputenc}
\usepackage[usenames,dvipsnames,svgnames,table]{xcolor}
\usepackage[all]{xy}
\usepackage{graphicx}
\usepackage{mathrsfs}
\usepackage{hyperref}

\def\be{\begin{equation}}
\def\ee{\end{equation}}
\def\bea{\begin{eqnarray}}
\def\eea{\end{eqnarray}}
\def\bal{\begin{align}}
\def\eal{\end{align}}

\newtheorem{thm}{Theorem}[section]

\newtheorem{prop}[thm]{Proposition}
\newtheorem{defn}{Definition}[section]

\newtheorem{exmp}{Example}[section]

\newtheorem{re}{Remark}[section]

\newcommand\TM{Teichm\"uller~}

\newcommand{\ud}{\mathrm{d}}
\newcommand{\ub}{\mathrm{b}}
\newcommand{\R}{\mathrm{R}}
\newcommand{\NS}{\mathrm{NS}}

\newcommand{\vm}{{v^{(-)}}}
\newcommand{\vp}{{v^{(+)}}}

\newcommand{\ka}{\kappa}
\newcommand{\hka}{{\hat \kappa}}

\def\be{\begin{equation}}
\def\ee{\end{equation}}
\def\bea{\begin{eqnarray}}
\def\eea{\end{eqnarray}}
\def\bal{\begin{align}}
\def\eal{\end{align}}

\newcommand\bZ{\mathbb{Z}}

\newcommand\cA{\mathcal{A}}

\newcommand{\hH}{{\hat H}}

\def\be{\begin{equation}}
\def\ee{\end{equation}}
\def\bea{\begin{eqnarray}}
\def\eea{\end{eqnarray}}
\def\bal{\begin{align}}
\def\eal{\end{align}}

\definecolor{lightblue}{rgb}{0.06, 0.75, 0.99}
\definecolor{darkblue}{rgb}{0.06, 0.3, 0.57}
\definecolor{darkgreen}{rgb}{0.0, 0.5, 0.0}
\definecolor{brightpink}{rgb}{1.0, 0.0, 0.5}
\definecolor{darkgreen}{rgb}{0.11, 0.35, 0.02}

\title{Heisenberg double and Drinfeld double of the quantum superplane}

\author[a,b]{Nezhla Aghaei}
\author[c]{Michal Pawelkiewicz} 

\affiliation[a]{Max Planck Institut f\"ur Mathematik, Vivatsgasse 7, 53111 Bonn, Germany.}
\affiliation[b] {Albert Einstein Center for Fundamental Physics, Institute for Theoretical Physics, University of Bern, Sidlerstrasse 5, Bern, ch-3012, Switzerland.}

\affiliation[c]{ Institut de Physique Theorique, CEA Saclay, 91191 Gif Sur Yvette, France.}

\emailAdd{nezhla.aghaee@mpim-bonn.mpg.de,~michal.pawelkiewicz@ipht.fr}

\abstract{
We study infinite dimensional generalisations of the Heisenberg doubles of the Borel half of $U_q(sl(2))$ and of $U_q(osp(1|2))$ and find associated canonical elements which satisfy pentagon equation. The former reproduces the canonical element, expressed using the Faddeev's quantum dilogarithm, which has been found by Kashaev to be realised within quantised \TM theory, while for the latter we show that it corresponds to an operator from quantised super \TM theory. We study infinite dimensional representations of those two Heisenberg doubles and, using an algebra homomorphism between Heisenberg doubles and Drinfeld doubles, we find associated representations of Drinfeld doubles of the Borel half of $U_q(sl(2))$ and of $U_q(osp(1|2))$. Moreover, we reproduce the previously obtained $R$-matrix for the former and derive a novel $R$-matrix for the latter representation.}

\begin{document} 

\maketitle
\flushbottom

\section{Introduction}

The methods comming from the representation theory of quantum groups have found a wide range of applications to the mathematical and theoretical physics. Quantum groups are relevant in conformal field theory \cite{Witten:1989rw}, where the algebras of screening charges and vertex operators satisfy the relations of q-deformed Lie algebras, which are themselves a well studied family of quantum groups. The fusion matrices of the conformal field theories were realised as a 6j symbols for representations of the associated quantum groups. Moreover, in the context of quantum integrable systems a systematic method of obtaining scatterring matrices \cite{zam} has been developed. These systems satisfied the so-called \textit{Yang-Baxter equation}~\cite{Baxter:1972hz, Yang:1967bm}
\begin{align}
R_{12}R_{13}R_{23} = R_{23}R_{13}R_{12},
\end{align}
There exist a systematic procedure to obtain solutions to the Yang-Baxter equation, which is based on the so-called \textit{Drinfeld double construction} \cite{DD,Majid,Majid0,Majidbook}. It allows to associate a new, quasi-triangular Hopf algebra, i.e. a Hopf algebra that admits a universal $R$-matrix, which satisfies the Yang-Baxter equation, to an arbitrary Hopf algebra and its dual. 

Another existing double construction, called the \textit{Heisenberg double construction} \cite{lu,Ka3}, admits a canonical element that satisfies not the Yang-Baxter equation, but rather a \textit{pentagon equation}~\cite{Maillet}
\begin{align}
S_{12}S_{13}S_{23}=S_{23}S_{12}.
\end{align}
Using Heisenberg doubles one can design the representations of Drinfeld doubles, as one can embed the elements of the Drinfeld double into a tensor square of Heisenberg doubles \cite{lu, Ka3}.

In mathematical physics, Heisenberg doubles appeared in particular in the context of the quantum \TM theory of Riemann surfaces \cite{Teschner:2005bz, Frenkel:2010bv, Ka1,Ka2,Ka4, Chekhov:1999tn}. The \TM theory is the theory of the deformations of the complex structures on Riemann surfaces. On the space of complex structures one can define the local coordinates using the triangulations of Riemann surfaces --- however, descriptions given by equivalent triangulations should be related by similarity transformations. The transfromation that realises this, i.e. the so-called flip move, relates different triangulations of a quadrilateral and is one of the generators of the Ptolemy grupoid.  The canonical element of the Heisenberg double of the quantum plane (that is, of a Borel half of $U_q(sl(2))$) evaluated on a particular family of infinite dimensional representations realises this flip move \cite{Ka1}. Moreover, the Fock coordinates associated to the edges of a quadrilateral correspond to the elements of the aformentioned Heisenberg double. 

From the Heisenberg double of the quantum plane Kashaev obtained a class of representations of the $U_q(sl(2))$ quantum group, as well as its associated $R$-matrix. Those representations were identified to be the class of infinite dimensional representations $\mathcal{P}_\alpha$ studied \cite{PT1,PT2} in connection with the Liouville field theory, which constitutes a prototypical non-trivial example of the noncompact conformal field theory \cite{T03, Teschner:2003em}. Using the means of harmonic analysis, Ponsot and Teschner investigated their properties. They have shown that the relation between the fusion category of the conformal field theory and the representation category of quantum group holds in the case of the Liouville field theory and $U_q(sl(2))$ quantum group. Moreover, the consistency of the bootstrap for the Liouville theory, i.e. the fact that the crossing-symmetry equation is satisfied by the three point function, was verified \cite{PT2,PT1,teschner}.

The canonical element $S$ in the \TM theory context was expressed in terms of a particular hyperbolic special function called the \textit{Faddeev's quantum dilogarithm} \cite{Faddeev:1993rs,Y.:2005tv},
\begin{equation} \label{dilog}
e_\ub(z) = \exp\left( \int_C \frac{e^{-2izw}}
{\sinh(w\ub)\sinh(w/\ub)}
\frac{\ud w}{4w}\right) \ ,
\end{equation}
and it can be regarded as a quantisation of the Roger's dilogarithm. In fact, the quantum dilogarithm has many elegant properties. In particular, the pentagon equation that it satisfies
\begin{equation}
e_\ub(\mathsf{p})e_\ub(\mathsf{x}) = e_\ub(\mathsf{x})e_\ub(\mathsf{p}+\mathsf{x})e_\ub(\mathsf{p}),
\end{equation}
for non-commutative variables $\mathsf{x},\mathsf{p}$ such that $[\mathsf{p},\mathsf{x}] = \frac{1}{2\pi i}$, is a generalisation of the five-term relation for Roger's five-term identity. The quantum dilogarithm also found applications in conformal field theory, topological field theory and hyperbolic geometry.

In the context of super \TM theory \cite{APT}, one can define the supersymmetric analogues of the Faddeev's quantum dilogarithm function 
\begin{align}\begin{aligned}
e_\R(x) &= e_\ub\left(\frac{x}{2}+\frac{i }{4}(\ub-\ub^{-1})\right)e_\ub\left(\frac{x}{2}-\frac{i }{4}(\ub-\ub^{-1})\right),\\
e_\NS(x) &= e_\ub\left(\frac{x}{2}+\frac{i }{4}(\ub+\ub^{-1})\right)e_\ub\left(\frac{x}{2}-\frac{i }{4}(\ub+\ub^{-1})\right) .
\end{aligned}\end{align}
For self-adjoint operators $\mathsf{p},\mathsf{x}$ such that $[\mathsf{p},\mathsf{x}]=\frac{1}{ \pi i}$ those supersymmetric quantum dilogarithms have been shown to satisfy four pentagon relations
\begin{subequations}\label{superpentagon}
	\begin{align}
	f_+(\mathsf{p}) f_+(\mathsf{x}) &=  f_+(\mathsf{x})f_+(\mathsf{x}+\mathsf{p})f_+(\mathsf{p}) -i f_-(\mathsf{x})f_-(\mathsf{x}+\mathsf{p})f_-(\mathsf{p}), \\
	f_+(\mathsf{p}) f_-(\mathsf{x}) &=  -i f_+(\mathsf{x})f_-(\mathsf{x}+\mathsf{p})f_-(\mathsf{p}) + f_-(\mathsf{x})f_+(\mathsf{x}+\mathsf{p})f_+(\mathsf{p}), \\
	f_-(\mathsf{p}) f_+(\mathsf{x}) &=  f_+(\mathsf{x}) f_+(\mathsf{x}+\mathsf{p}) f_-(\mathsf{p}) -i f_-(\mathsf{x})f_-(\mathsf{x}+\mathsf{p})f_+(\mathsf{p}) , \\
	f_-(\mathsf{p}) f_-(\mathsf{x}) &=  i f_+(\mathsf{x}) f_-(\mathsf{x}+\mathsf{p}) f_+(\mathsf{p}) - f_-(\mathsf{x})f_+(\mathsf{x}+\mathsf{p})f_-(\mathsf{p}),
	\end{align}
\end{subequations}
where $f_\pm(x) = e_\R(x) \pm e_\NS(x)$. As a consequence, the supersymmetric analogue of the flip operator satisfies the (graded) pentagon equation. \\

In this paper, we aim at constructing the canonical elements which satisfy pentagon equation in terms of the elements of continuous Heisenberg doubles of the Borel half of $U_q(sl(2))$ and $U_q(osp(1|2))$. The fact that the algebras considered are spanned by a continuum of basis elements require form us a careful treatment when compared to the construction for the discrete cases. As a result, we are able to obtain the canonical elements in terms of the basis elements of Heisenberg algebras and, using representations, to arrive at formulae which were found to be present in the \TM and super \TM theories. Moreover, using those representations we construct the representations of the Drinfeld doubles, as well as realise the associated $R$-matrices. The expression for the $U_q(osp(1|2))$ $R$-matrix is a new one. \\

The paper is organised as follows. In the section \ref{chapter2} we present a general theory of ($\mathbb{Z}_2$-graded) Heisenberg doubles and consider some discrete algebras as demonstrative examples. In the section \ref{chapter3} we define the infinite dimensional, continuous version of the Borel half of $U_q(sl(2))$ and construct its Heisenberg double along with the canonical element corresponding to the one obtained previously by Kashaev \cite{Ka3}. We also consider an infinite dimensional representation of this algebra, which was found to be relevant for \TM theory. In the section \ref{chapter4} we consider the continuous analogue of the Borel half of $U_q(osp(1|2))$ and we construct its Heisenberg double. We stress how the construction differs from the one in the previous section as a consequence of the non-trivial grading of $U_q(osp(1|2))$. We find the canonical element and represent it in a way which appeared previously in super \TM theory. In the section \ref{chapter5} we quickly recall some information about the Drinfeld double construction. We generalise the algebra homomorphism between Heisenberg and Drinfeld doubles, given in \cite{Ka3}, to the $\mathbb{Z}_2$-graded case, and consider the representations of the continuous versions of the Drinfeld doubles which stem from the Heisenberg double representations. We derive also previously obtained $R$-matrix in the $U_q(sl(2))$ case and a new $R$-matrix for the $U_q(osp(1|2))$ one. Section \ref{chapter6} offers an outlook.
\

\section{Heisenberg doubles}\label{chapter2}

In this chapter we will shortly describe the basic notions about $\mathbb{Z}_2$-graded Heisenberg doubles. We will sketch how the Heisenberg double construction works. The exposition is structured in a way that is similar to \cite{Ka3}, it is however generalised to work in the $\mathbb{Z}_2$-graded setting.

We start from a short description of $\mathbb{Z}_2$-graded Hopf algebras, and using the Hopf action we will define a smash product of a Hopf algebra with its dual. 
At the end of this chapter we will illustrate the Heisenberg double construction  on the examples of discrete versions of $U_q(sl(2))$ and $U_q(osp(1|2))$ Heisenberg doubles. They constitute discrete prototypes for the continuous Heisenberg doubles of the Borel halfs of $U_q(sl(2))$ and $U_q(osp(1|2))$, which will be considered in chapters \ref{chapter3} and \ref{chapter4} respectively.\\

Let us consider a $\mathbb{Z}_2$-graded Hopf-algebra $(\mathcal{A},m,\eta,\Delta,\epsilon,\gamma)$, where $\mathcal{A}$ is a $\mathbb{Z}_2$-graded vector space equipped with the multiplication $m:\mathcal{A}\otimes\mathcal{A}\to\mathcal{A}$, an unit $\eta:\mathbb{C}\to \cA$, a co-multiplication $\Delta: \cA\to\cA\otimes\cA$, a co-unit $\epsilon:\cA\to\mathbb{C}$ and an antipode $\gamma:\cA\to\cA$. $\cA$ decomposes into a direct sum of two sub-spaces $\cA = \cA_0 \oplus \cA_1$, which are called even and odd respectively. We denote the degree of an element $x\in\cA_i$ by $|x|=i$, and we will call an element $x$ even if $|x|=0$ and odd otherwise. The graded tensor product of two algebras $\cA$ and $\mathcal{B}$ is then defined by the following equation for $a_1, a_2\in\cA$, $b_1,b_2\in \mathcal{B}$,
\begin{align}
	(a_1 \otimes b_1)\cdot(a_2 \otimes b_2) = (-1)^{|b_1||a_2|} a_1 a_2 \otimes b_1 b_2.
\end{align}
 The maps $m$ and $\eta$ are subjected to the associativity and unitality relations
\begin{align}\label{hopf_alg_axioms_1}
	m \circ (m\otimes id) &= m \circ (id \otimes m),\\
	m \circ (\eta \otimes id) &= id = m \circ (id \otimes \eta) ,
\end{align}
while maps $\Delta$ and $\epsilon$ have to satisfy the co-associativity and co-unitality relations
\begin{align}
	(\Delta\otimes id)\circ \Delta &= \Delta\circ (id \otimes \Delta),\\
	(\epsilon \otimes id)\circ \Delta &= id = (id \otimes \epsilon)\circ\Delta.
\end{align}
Moreover, the co-product $\Delta$ and co-unit $\epsilon$ are algebra homomorphisms, and the antipode $\gamma$ is a graded algebra anti-homomorphism and a graded co-algebra anti-homomorphism which satisfies the relations
\begin{align}\label{hopf_alg_axioms_3}
	m \circ (id\otimes \gamma)\circ \Delta = m\circ (\gamma\otimes id)\circ \Delta = \eta\circ\epsilon .
\end{align}
All the above maps are grade preserving.  \\

Moreover, we consider a Hopf algebra $(\mathcal{A}^*,{\hat m},{\hat \eta},{\hat \Delta},{\hat \epsilon},{\hat \gamma})$ which is dual to $\cA$. The Hopf algebras $\mathcal{A}$ and $\mathcal{A}^*$ are dual in a sense that the vector spaces $\mathcal{A}$ and $\mathcal{A}^*$ are dual as vector spaces, and there exists a non-degenerate duality pairing (also called a \textit{Hopf pairing}) $(,): \mathcal{A}\times\mathcal{A}^*\rightarrow \mathbb{C}$, for which the following relations are satisfied
\begin{gather}\label{duality_bracket_eq1}
 (x, f g) = (\Delta(x),f\otimes g), \qquad (x y, f) = (x\otimes y, {\hat \Delta}(f)), 
\end{gather}
between multiplications and co-multiplications,
\begin{gather}
 (\eta(1) , f) = {\hat \epsilon}(f), \qquad \epsilon(x) = (x,{\hat \eta}(1)), 
\end{gather}
between unit and co-unit maps,
\begin{gather}\label{duality_bracket_eq3}
 (\gamma(x),f) = (x,{\hat \gamma}(f)) ,
\end{gather}
and between antipodes, where
\begin{align*}
 (x\otimes y, f\otimes g) = (-1)^{|y||f|} (x,f)(y,g) ,
\end{align*}
for $x,y\in\mathcal{A}$, $f,g\in\mathcal{A}^*$. \\

The ordinary tensor product $\cA^* \otimes \cA$ has a straight-forward product given by $(1\otimes x)(f\otimes 1) = x\otimes f $. However, in order to construct a Heisenberg double, we are interested in equipping the space $\cA^* \otimes \cA$ with a non-trivial algebra structure between the elements which belong to the subalgebras $\cA$ and $\cA^*$. In order to achieve that, we will use the Hopf pairing $(,)$ to define a left action of $\cA$ on $\cA^*$ and consequently a smash product algebra $\mathcal{A}^* \rtimes \mathcal{A}$. \\

Using the duality pairing $(,)$ we can define a left action $\triangleright$ of a Hopf-algebra $\mathcal{A}$ on $\mathcal{A}^*$ given by
\begin{equation}\label{action_of_cA_on_cA*}
 x\otimes f \mapsto \sum_{(f)} (-1)^{|f_{(1)}||f_{(2)}|} (x, f_{(2)}) f_{(1)} =: x \triangleright f ,
\end{equation}
where $x\in\cA, f\in\cA^*$ and where we denote the coproduct ${\hat \Delta}(f) = \sum_{(f)} f_{(1)}\otimes f_{(2)}$ using usual Sweedler notation \cite{Sw}. The action~\eqref{action_of_cA_on_cA*} makes $\cA^*$ into a module algebra over the Hopf algebra $\cA$, i.e.\
the action is compatible with the multiplication in $\cA^*$ in the sense that
\begin{equation}\label{eq:x-fg}
 x \triangleright (f g) = \sum_{(x)} (-1)^{|f| |x_{(2)}|} (x_{(1)} \triangleright f) \, (x_{(2)} \triangleright g) ,
\end{equation}
where $x\in\cA$ and $f,g\in\cA^*$. Using a left action $\triangleright$ one can construct a smash product algebra $H(\cA) = \cA^* \rtimes \cA$ by defining the multiplication
\begin{equation}\label{heisenbergmultiplication}
(f\otimes x)(g\otimes y) = \sum_{(x)} (-1)^{|g||x_{(2)}|} f (x_{(1)} \triangleright g)\otimes x_{(2)} y ,
\end{equation}
where $x,y\in \mathcal{A}$, $f,g\in \mathcal{A}^*$. 

\begin{defn}
A {\normalfont{Heisenberg double}} of a Hopf algebra $\cA$ is the smash product algebra $H(\cA) = \cA^* \rtimes \cA$ with the multiplication given by the equation \eqref{heisenbergmultiplication}.
\end{defn}

The Heisenberg double has $\mathcal{A}$ and $\mathcal{A}^*$ as subalgebras through canonical embeddings $\mathcal{A}^*\ni f\mapsto f\otimes1\in H(\mathcal{A})$ and $\mathcal{A}\ni x \mapsto 1\otimes x \in H(\mathcal{A})$. \\

It will be convienient to recast the definition above in a basis dependent way. In order to do that, we first choose a basis of $\cA$. The basis will be given by a collection of vectors $\{e_\alpha\}_{\alpha\in I}$, where $I$ is a (possibly infinite) set. Then, the multiplication and co-multiplication of the basis elements is given by
\begin{align}
& e_\alpha e_\beta = \sum_{\gamma\in I} m^\gamma_{\alpha\beta} e_\gamma,
&& \Delta(e_\alpha) = \sum_{\beta, \gamma\in I} \mu^{\beta\gamma}_\alpha e_\beta\otimes e_\gamma ,
\end{align}
where $m^\gamma_{\alpha\beta}$ and $\mu^{\beta\gamma}_\alpha$ are respectively multiplication and comultiplication coefficients. With the choice of a basis $\{e^\alpha\}_{\alpha\in I}$ of $\cA^*$ dual to $\{e_\alpha\}_{\alpha\in I}$ in the sense
\begin{equation}\label{duality}
(e_\alpha,e^\beta) = \delta_\alpha^\beta ,
\end{equation}
the multiplication and co-multiplication on $\cA^*$ are as follows
\begin{align}
&e^\alpha e^\beta = \sum_{\gamma\in I} (-1)^{|\alpha||\beta|} \mu_\gamma^{\alpha\beta} e^\gamma,
&&{\hat \Delta}(e^\alpha) = \sum_{\beta, \gamma\in I} (-1)^{|\beta||\gamma|} m_{\beta\gamma}^\alpha e^\beta\otimes e^\gamma .
\end{align}
The Heisenberg double $H(\mathcal{A})$ is thus spanned by basis elements $\{ e^\alpha \otimes e_\beta\}_{\alpha,\beta\in I}$, written in terms of the basis elements of $\cA$ and $\cA^*$. The multiplication \eqref{heisenbergmultiplication} on basis elements $e^\alpha \otimes e_\beta$ has the following form
\begin{equation}
 (e^\alpha \otimes e_\beta)(e^\gamma \otimes e_\delta) = \sum_{\epsilon,\pi,\rho,\sigma,\tau\in I} (-1)^{|\beta||\gamma| + |\pi||\epsilon| + |\pi||\alpha| + |\epsilon|} m^\gamma_{\pi\epsilon} m^\tau_{\rho\delta} \mu^{\epsilon\rho}_\beta \mu^{\alpha\pi}_\sigma e^\sigma \otimes e_\tau.
\end{equation} 

It is important to note that the Heisenberg double $H(\cA)$ is not a Hopf algebra. The algebra structure given by $\eqref{heisenbergmultiplication}$ is not compatible with the co-products on $\Delta$, ${\hat \Delta}$ defined on the initial Hopf algebras $\cA$ and $\cA^*$. In this it differs from the Drinfeld double algebra (which will be discussed in section \ref{chapter5}), which is a (quasi-triangular) Hopf algebra, and not only an algebra.

\begin{defn}
For Heisenberg algebra of interest to us is a \normalfont{canonical element} $S \in H(\mathcal{A})\otimes H(\mathcal{A})$
\begin{equation}\label{canonicalelementdef}
S = \sum_{\alpha\in I} (-1)^{|\alpha|} (1\otimes e_\alpha)\otimes (e^\alpha \otimes 1) ,
\end{equation} 
\end{defn}

\begin{prop}
The canonical element $S$ satisfies the \normalfont{graded pentagon relation}
\begin{equation}\label{pentagonrelation} 
S_{12}S_{13}S_{23} = S_{23}S_{12},
\end{equation}
where we use a notation for which $S_{12} = S \otimes (1\otimes1)$, $S_{23} = (1\otimes1)\otimes S$ and $S_{13} = \sum_{\alpha\in I} (-1)^{|\alpha|} e^\alpha \otimes (1\otimes1) \otimes e_\alpha$. 
\end{prop}

Although, as we have mentioned previously, the Heisenberg double is not a Hopf algebra, the canonical element $S$ does encode the co-products on the initial Hopf algebras $\mathcal{A}$ and $\mathcal{A}^*$ in the following way
\begin{align}\begin{aligned}
\Delta(e_\alpha) &= S^{-1} \big( (1\otimes1)\otimes (1\otimes e_\alpha) \big) S , \\
{\hat \Delta}(e^\alpha) &= S\big((e^\alpha \otimes 1) \otimes (1\otimes1)\big)S^{-1} ,
\end{aligned}\end{align}
where $\Delta(e_\alpha)$ and ${\hat \Delta}(e^\alpha)$ should be understood as elements of the Heisenberg double obtained by embedding the co-products under the canonical embeddings $\cA\hookrightarrow H(\cA)$ and $\cA^*\hookrightarrow H(\cA)$. \\

\begin{re}
 To keep the notation compact, from now on we will denote the elements $1\otimes e_\alpha $ and $e^\alpha\otimes1 $ of the Heisenberg double $H(\cA)$ simply as $e_\alpha$ and $e^\alpha$ respectively.
\end{re}

\begin{exmp} Heisenberg double of the Borel half of $U_q(sl(2))$ \end{exmp} \vspace{-0.25cm}

As an instructive example, lets us consider a Heisenberg double of the Borel half of $U_q(sl(2))$, which should be considered as a discrete prototype of the continuous algebra considered in chapter~\ref{chapter3}. The Borel half algebra $\mathcal{A}=\mathcal{B}(U_q(sl(2))) = \{ H^m E^n \}_{m,n=0}^\infty$ is generated by elements $H, E$ with a commutation relation
\begin{align}\label{discrete_multiplication}
 &[H,E] = -ib E,
 \end{align}
and a coproduct as follows
 \begin{align}\label{coproductsltwo}
 &\Delta(H) = H\otimes1 + 1\otimes H, && \Delta(E) = E\otimes e^{2\pi b H} + 1\otimes E,
\end{align}
where $q = e^{\pi i b^2}$ is the deformation parameter. In addition, the antipode is
\begin{align}
\begin{aligned}
 \gamma(H) &= -H, && \gamma(E) &= q E e^{-2\pi b H}.
 \end{aligned}
\end{align}
We can choose the basis elements of $\mathcal{A}$ in the following way
\begin{equation}
 e_{m,n} = \frac{q^{n}}{m! (q^2)_n} (i b^{-1} H)^m (i E)^n,
\end{equation}
where q-numbers $(q)_n$ are defined as $(q)_n = (1-q)\ldots(1-q^n)$ and $n,m\in\mathbb{N}$.

Using the properties of the generators $H, E$ and the binomial and q-binomial formulae one can find the multiplication and co-multiplication of the basis elements
\begin{align*}
 & e_{m,n} e_{k,l} = \sum_{j=0}^k {m+j\choose j} {n+l\choose l}_{q^2} \frac{(-n)^{k-j}}{(k-j)!} e_{m+j,n+l}, \\
 &\Delta(e_{m,n}) = \sum_{k=0}^m \sum_{l=0}^n \sum_{p=0}^\infty {k+p\choose k} (n-l)^p (-2\pi i b^2)^p e_{m-k,n-l}\otimes e_{k+p,l}.
\end{align*}
where ${n\choose k}$ is an ordinary and ${n\choose k}_q = \frac{(q)_{n}}{(q)_{k}(q)_{n-k}}$ is q-deformed binomial coefficient.
 
The dual Borel half algebra $\mathcal{A}^*=\mathcal{B}(U_q(sl(2)))=\{{\hat H}^m F^n\}_{m,n=0}^\infty$ on the other hand is generated by elements $\hat H, F$ which satisfy a commutation relation
\begin{align}\label{discrete_multiplication_dual}
 &[\hat H,F] = +ib F,
 \end{align}
 and have a coproduct given by
 \begin{align}\label{discrete_comultiplication_dual}
 &{\hat \Delta}(\hat H) = \hat H\otimes1 + 1\otimes \hat H,
 &&{\hat \Delta}(F) = F\otimes e^{-2\pi b {\hat H}} + 1\otimes F.
\end{align}
The antipode is
\begin{align}
\begin{aligned}
 {\hat\gamma}(\hat H) &= - \hat H, && {\hat\gamma}(F) &= q F e^{2\pi b \hat H}.
 \end{aligned}
\end{align}
The basis elements for $\mathcal{A}^*$ are given by
\begin{equation}
 e^{m,n} = (2 \pi b {\hat H})^m (i F)^n ,
\end{equation}
and their multiplication and co-multiplication is as follows
\begin{align*}
 & e^{m,n} e^{k,l} = \sum_{j=0}^k {k\choose j} (n)^{k-j} (-2\pi i b^2)^{k-j} e^{m+j,n+l} , \\
 &{\hat \Delta}(e^{m,n}) = \sum_{k=0}^m \sum_{l=0}^n \sum_{p=0}^\infty {m\choose k} {n\choose l}_{q^2} \frac{(-n+l)^p}{p!} e^{m-k,n-l}\otimes e^{k+p,l}.
\end{align*}
By inspection the multiplication and co-multiplication coefficients are equal to
\begin{align*}
 & m^{r,s}_{m,n;k,l} = {r\choose r-m} {n+l\choose l}_{q^2} \frac{(-n)^{k-r+m}}{(k-r+m)!} \Theta(r-m) \Theta(k-r+m) \delta_{s,n+l}, \\
 & \mu_{r,s}^{m,n;k,l} = {k\choose r-m} (n)^{k-r+m} (-2\pi i b^2)^{k-r+m} \Theta(r-m) \Theta(k-r+m) \delta_{s,n+l} ,
\end{align*}
and that $\cA$ and $\cA^*$ are dual as Hopf algebras. The Heisenberg double $H(\cA)$ is then given by the generators $H,\hat H, E, F$, which satisfy the commutation relations \eqref{heisenbergmultiplication}, which are as follows
\begin{align}\label{commutationrelation}
\begin{aligned}
 &[H,E] = -ibE, &&\qquad [H,F] = +ibF,\\
 &[\hat H, E] = 0, &&\qquad [\hat H,F] = +ibF,\\
 &[H,\hat H] = \frac{1}{2\pi i}, &&\qquad [E,F] = (q-q^{-1})e^{2\pi b H} .
 \end{aligned}
\end{align}
The canonical element $S$ is given in terms of generators as
\begin{align}
 &S = \exp(2 \pi i H\otimes\hat H) (- q E\otimes F;q^2)^{-1}_\infty ,
\end{align}
where the special function $(x;q)_\infty$, known under a name of a \textit{quantum dilogarithm}, is defined as
\begin{align}\label{specialfunction}
 &(x;q)_\infty^{-1} = \prod_{k=0}^\infty \frac{1}{1-xq^k}=  \sum_{k=0}^\infty \frac{x^k}{(q)_k} .
\end{align}
Using the properties of the quantum dilogarithm function one can check explicitly that the pentagon equation is satisfied. In particular, it reduces to the identity
\begin{align}\label{identity}
 (V;q^2)_\infty (U;q^2)_\infty = (U;q^2)_\infty ((1-q^2)^{-1} [U,V];q^2)_\infty (V;q^2)_\infty ,
\end{align}
for $U = -q1\otimes E\otimes F,V=-q E\otimes F\otimes 1$, which was shown to be satisfied by the quantum dilogarithm $(x;q)_\infty$.
The square brackets denote the commutator, and operators $U$ and $V$ satisfy the following algebraic relations \cite{Ka3}
\begin{align}
&W= UV-q^2VU,&&[U,W] = [V,W] = 0.
\end{align}

\newpage
\begin{exmp} Heisenberg double of the Borel half of $U_q(osp(1|2))$ \end{exmp} \vspace{-0.25cm}
 
As another informative example, lets us consider a Heisenberg double of the Borel half of $U_q(osp(1|2))$, which should be considered as a discrete prototype of the continuous algebra considered in chapter~\ref{chapter4}. 
The Borel half algebra $\mathcal{A}=\mathcal{B}(osp(1|2)))) = \{ H^m \vp^n \}_{m,n=0}^\infty$ is generated by an even graded element $H$ and an odd graded element $\vp$ with a commutation relation
\begin{align}
	&[H,\vp] =-ib \vp,
\end{align}
a coproduct 
\begin{align}\label{coproductsltwo}
	&\Delta(H) = H\otimes1 + 1\otimes H, &&  \Delta(\vp) = \vp\otimes e^{\pi b H} + 1\otimes \vp, 
\end{align}
and the antipode as follows
\begin{align}
\begin{aligned}
 \gamma(H) &= -H, && \gamma(\vp) = q^\frac{1}{2} \vp e^{-\pi b H} ,
 \end{aligned}
\end{align}
where $q = e^{i\pi b^2}$ is the deformation parameter. We can choose the basis elements of $\mathcal{A}$ in the following way
\begin{equation}
	e_{m,n} = (-1)^{n(n-1)/2} \frac{q^{n/2}}{m! (-q)_n} (i b^{-1} H)^m (i \vp)^n .
\end{equation}
Using the properties of the generators $H, \vp$ and the binomial and q-binomial formulae one can find the multiplication and co-multiplication of the basis elements
\begin{align}
 & e_{m,n} e_{k,l} = \sum_{j=0}^k {m+j\choose j} {n+l\choose l}_{-q} \frac{(-n)^{k-j}}{(k-j)!} e_{m+j,n+l},\label{discrete_osp_multiplication_basis} \\
 &\Delta(e_{m,n}) = \sum_{k=0}^m \sum_{l=0}^n \sum_{p=0}^\infty {k+p\choose k} (n-l)^p (-\pi i b^2)^p e_{m-k,n-l}\otimes e_{k+p,l}. \label{discrete_osp_comultiplication_basis}  
\end{align}
The dual Borel half algebra $\mathcal{A}^*=\mathcal{B}(U_q(osp(1|2)))=\{{\hat H}^m \vm^n\}_{m,n=0}^\infty$ on the other hand is generated by elements $\hat H, \vm$ which satisfy a commutation relation
\begin{align}
	[\hat H,\vm] = +ib\vm,
\end{align}
and have a coproduct given by
\begin{align}
	&{\hat \Delta}(\hat H) = \hat H\otimes1 + 1\otimes \hat H,
	&&{\hat \Delta}(\vm) = \vm\otimes e^{-\pi b {\hat H}} + 1\otimes \vm,
\end{align}
with the antipode
\begin{align}
\begin{aligned}
 \hat\gamma(\hat H) &= - \hat H, && \hat\gamma(\vm) = q^{-\frac{1}{2}} \vm e^{\pi b \hat H}. 
 \end{aligned}
\end{align}
The basis elements for $\mathcal{A}^*$ are given by
\begin{equation}
 e^{m,n} = (\pi b {\hat H})^m (i v^{(-)})^n , 
\end{equation}
and their multiplication and co-multiplication are as follows
\begin{align}
 & e^{m,n} e^{k,l} = \sum_{j=0}^k {k\choose j} (n)^{k-j} (-\pi i b^2)^{k-j} e^{m+j,n+l} , \label{discrete_osp_multiplication_basis_dual} \\
 &{\hat \Delta}(e^{m,n}) = \sum_{k=0}^m \sum_{l=0}^n \sum_{p=0}^\infty {m\choose k} {n\choose l}_{-q} \frac{(-n+l)^p}{p!} e^{m-k,n-l}\otimes e^{k+p,l}. \label{discrete_osp_comultiplication_basis_dual} 
\end{align}
By inspection the multiplication and co-multiplication coefficients are equal to
\begin{align*}
 & m^{r,s}_{m,n;k,l} = {r\choose r-m} {n+l\choose l}_{-q} \frac{(-n)^{k-r+m}}{(k-r+m)!} \Theta(r-m) \Theta(k-r+m) \delta_{s,n+l}, \\
 & \mu_{r,s}^{m,n;k,l} = {k\choose r-m} (n)^{k-r+m} (-\pi i b^2)^{k-r+m} \Theta(r-m) \Theta(k-r+m) \delta_{s,n+l} , 
\end{align*}
and that $\cA$ and $\cA^*$ are dual as Hopf algebras. 
The discrete Heisenberg double can be defined as an algebra generated by the even elements $H$ and $\hat H$ and the odd elements $\vp$ and $\vm$ satisfying (anti-)commutation relations
\begin{align}\label{discrete_osp_heisenberg_commutation_relations}\begin{aligned}
		& [H,\hat H] = \frac{1}{\pi i}, && \{\vp, \vm\} = -e^{\pi bH}(q^{\frac{1}{2}}+q^{-\frac{1}{2}}) ,\\
		& [H,\vp]=-ib\vp, && [H, \vm] = ib \vm, \\
		& [\hat H, \vp] = 0, && [\hat H,\vm]=+ib \vm.
\end{aligned}\end{align}
The canonical element $S$ is given in terms of generators as
\begin{align}
	&S = \exp(\pi i H\otimes {\hat H}) (-q^\frac{1}{2} \vp\otimes \vm;-q)^{-1}_\infty. 
\end{align}
Using the properties of the quantum dilogarithm function one can check explicitly that the pentagon equation is satisfied. In particular, it reduces to the identity \eqref{identity} for $U = -q^\frac{1}{2} 1 \otimes \vp \otimes \vm,V=-q^\frac{1}{2} \vp\otimes \vm \otimes 1$ and with the replacement $q^2\to -q$.  The operators $U$ and $V$ satisfy the following algebraic relations
\begin{align}
&W= UV+qVU,&&[U,W] = [V,W] = 0.
\end{align}

\section{Continuous Heisenberg double of the Borel half of $U_q(sl(2))$}\label{chapter3}
In this section we will provide a discussion of a continuous version of the Heisenberg double of the Borel half of $U_q(sl(2))$. We will describe the multiplication and co-multiplication structures of the continuous Borel half, which follow from the commutation relations and the co-product for the generators of the discrete algebra and the continuous version of the binomial and q-binomial formulae. Afterwards, we construct the canonical element $S$ satisfying the pentagon equation, which is expressed using the Faddeev's quantum dilogarithm. At the end of the section we recall a particular representation of continuous Heisenberg double which was introduced by Kashaev \cite{Ka2}.\\

As already described in the section \ref{chapter2}, the discrete version of the Heisenberg double $H(\mathcal{A})$ of the Borel half $\mathcal{A} = \mathcal{B}(U_q(sl(2)))$ can be defined as an algebra spanned by the elements $\{{\hat H}^m F^n H^k E^l \}_{m,n,k,l=0}^\infty$ subjected to the commutation relations \eqref{commutationrelation}, where $q=e^{i\pi b^2}$ for a parameter $b$ such that $b^2\in\mathbb{R}\slash\mathbb{Q}$. 

In the case of the continuous version of the algebra, instead of integer powers of the generators, we are interested in the generators taken to the pure imaginary powers --- i.e. instead of $H^{i_1} E^{i_2}$ for $i_j\in\mathbb{Z}_{\geq0}$ we would like to consider $ H^{i \alpha_1} E^{i \alpha_2}$ for $\alpha_j\in\mathbb{R}$. This modification would not involve any subtleties if all the generators were positive (or at least non-negative) definite. The situation would be similar to taking a complex power of a positive real number, which does not require specifting the branch of the logarithm --- however, taking a complex power of a negative real number does.

Unfortunately, for the algebra $\cA$ not all generators are be positive. We assume that the generator $E$ will be realised as positive operator while the generator $H$ will not be. This positivity issue will be resolved by using principal value prescription for those generators which belong to the Cartan subalgebra of the Borel half and its dual. The principal value is given by
\begin{align}\begin{aligned}
		{(\epsilon y)}_{pv}^{is}&=|y|^{is}\Theta(\epsilon y)+e^{-\pi s}|y|^{is}\Theta(-\epsilon y),
	\end{aligned}
\end{align}
where $y>0$ and $\epsilon=\pm1$ and where $\Theta(y)$ is a Heaviside theta function. The Hopf algebra composed only of the Cartan subalgebra constitutes an instructive toy model, which because of its simplicity clarifies the construction for the full $\mathcal{B}(U_q(sl(2)))$. We refer to the appendix \ref{appendixB} for its discussion. \\

Starting from the equations \eqref{discrete_multiplication}, \eqref{coproductsltwo} one can derive the multiplication and co-multiplication relations for the elements of the form $ (\pm H)^{is}_{pv} E^{ib^{-1}t} $ using the continuous binomial formulae \eqref{continuousbinomialformulae}, the continuous version of the q-binomial formula \eqref{continuous_q_binomialformula} as well as the Mellin transform of the exponential function \eqref{mellin_transform_exponential}. The result of that calculation is presented below as the relations satisfied by the basis elements. \\

Let us start with the continuous version of the Borel half $\cA$ of $U_q(sl(2))$. We define it as an algebra $\cA$ spanned by the elements $\{e(s,\epsilon,t)\}_{s,t\in \mathbb{R}, \epsilon=\pm1}$, which satisfy the following multiplication 
\begin{small}\begin{align}\begin{aligned}\label{continuous_multiplication}
			& e(s,\epsilon,t) e(s',\epsilon,t')= \int \frac{\ud \sigma}{(2\pi)^2} {{-i(s+s'-\sigma)}\choose{-i s}}_{\Gamma} \Gamma(-i\sigma) { {i(t+t')} \choose {it}}_\ub e^{-2\pi i tt'} \times \\
			&\times |2\pi t|^{i\sigma} [e^{\pi \sigma/2 }\Theta(-\epsilon t) +\Theta(\epsilon t)e^{-\pi \sigma/2 }] e(s+s'-\sigma,\epsilon,t+t') ,\\
			& e(s,\epsilon,t) e(s',-\epsilon,t') = \int \frac{\ud \sigma}{(2\pi)^2} {{-i(s+s'-\sigma)}\choose{-i s}}_{\Gamma} \Gamma(-i\sigma) { {i(t+t')} \choose {it}}_\ub e^{-2\pi i tt'} \times \\
			&\times |2\pi t|^{i\sigma} [e^{\pi \sigma /2}\Theta( \epsilon t) +e^{-\pi \sigma /2}\Theta(-\epsilon t)] \left[ \frac{e^{-\pi(s'-\sigma)}(1-e^{-2\pi s})}{1-e^{-2 \pi(s'+s-\sigma)}} e(s+s'-\sigma,\epsilon,t+t') + \right. \\
			& \left. + \frac{e^{-\pi s}(1-e^{-2\pi (s'-\sigma)})}{1-e^{-2 \pi(s'+s-\sigma)}} e(s+s'-\sigma,-\epsilon,t+t') \right] ,
\end{aligned}\end{align}\end{small}
and co-multiplication relations
\begin{small}\begin{align}\begin{aligned}\label{continuous_comultiplication}
			& \Delta(e(s,\epsilon,t)) = \int\frac{\ud\sigma \ud{\tilde \sigma}}{2\pi}\ud\tau { -i(s-\sigma+{\tilde\sigma}) \choose -i{\tilde \sigma}}_\Gamma |\tau|^{i\tilde\sigma} e(\sigma,\epsilon,\tau) \otimes \\
			&\otimes \left\{ \left( \Theta(\epsilon \tau) + \frac{e^{-\pi\tilde\sigma}(1-e^{-2\pi(s-\sigma)})}{1-e^{-2\pi(s-\sigma+\tilde\sigma)}} \Theta(-\epsilon\tau) \right) e(s-\sigma+\tilde\sigma,\epsilon,t-\tau) + \right. \\
			&\left. +\frac{e^{-\pi(s-\sigma)}(1-e^{-2\pi\tilde\sigma})}{1-e^{-2\pi(s-\sigma+\tilde\sigma)}} \Theta(-\epsilon\tau) e(s-\sigma+\tilde\sigma,-\epsilon,t-\tau) \right\} ,
\end{aligned}\end{align}\end{small}
where ${ s \choose \sigma }_\Gamma = \frac{\Gamma(\sigma)\Gamma(s-\sigma)}{\Gamma(s)}$ is a continuous version of the binomial coefficent and $ { t \choose \tau }_\ub = \frac{G_\ub(-\tau)G_\ub(-t+\tau)}{G_\ub(-t)}$ is a continuous version of the q-binomial coefficient, defined using a special function $G_\ub$ which is related to the Fadeev's quantum dilogarithm (c.f. appendix \ref{appendixA})
\begin{equation}
	e_\ub(x) = \frac{\exp[-\frac{i\pi}{4}-\frac{i\pi}{12}(b^2+b^{-2})]}{G_\ub(\frac{Q}{2}- i x)}.
\end{equation}
The elements $e(s,\epsilon,t)$ admit a presentation in terms of the generators $H,E$ (that is the elements which generated the discrete version of the algebra in section \ref{chapter2}) in the following way
\begin{align}\label{basis}
	e(s,\epsilon,t) &= f_{\epsilon}(s,t) (\epsilon 2\pi H)^{is}_{pv} E^{ib^{-1}t} ,
\end{align}
where
\begin{align*}
	f_{\epsilon}(s,t) &= \frac{1}{2\pi} \Gamma(-is) G_\ub^{-1}(Q+it) e^{-\pi t Q} e^{\pi s/2} .
\end{align*}
One can reproduce the discrete multiplication and co-multiplication relations by analytically continueing the relations \eqref{continuous_multiplication} and \eqref{continuous_comultiplication} to the values $s=-im$, $t=-ibn$, where $n,m\in\mathbb{Z}_{\geq0}$. \\

After describing the continuous version of the algebra $\cA$, we focus on the algebra $\cA^*$ dual to it. Let us define it as being spanned by the elements $\{ {\hat e}(s,\epsilon,t)\}_{s,t\in \mathbb{R}, \epsilon=\pm1}$, with the following product
\begin{small}\begin{align}\begin{aligned}
			& {\hat e}(s,\epsilon,t){\hat e}(s',\epsilon,t') = \int \frac{\ud \sigma}{2\pi} {{-is'}\choose{-i\sigma}}_{\Gamma} |t|^{i\sigma} \left[ \Theta(-\epsilon t) \frac{e^{-\pi \sigma}(1-e^{-2\pi (s'-\sigma)})} {1-e^{-2\pi s'}}+\Theta(\epsilon t) \right] \times \\
			& \qquad\qquad \times {\hat e}(s+s'-\sigma,\epsilon,t+t') , \\
			& {\hat e}(s,\epsilon,t){\hat e}(s',-\epsilon,t') = \int \frac{\ud \sigma}{2\pi} {{-is'}\choose{-i\sigma}}_{\Gamma} |t|^{i\sigma} \Theta(-\epsilon t) \frac{e^{-\pi (s'-\sigma)}(1-e^{-2\pi \sigma})}{1-e^{-2\pi s'}} {\hat e}(s+s'-\sigma,\epsilon,t+t') , 
	\end{aligned}\end{align}
\end{small}and coproduct
\begin{small}\begin{align}\begin{aligned}
			&\Delta({\hat e}(s,\epsilon,t)) = \int \frac{\ud\sigma}{2\pi} \frac{\ud {\tilde \sigma}}{2\pi} \ud\tau \Gamma(-i{\tilde \sigma}) {{-is}\choose{-i\sigma}}_{\Gamma} { it \choose i\tau }_\ub e^{-2\pi i\tau(t-\tau)} |2\pi\tau|^{i{\tilde \sigma}} \times \\
			&\times \left\{ (e^{\pi{\tilde \sigma}/2} \Theta(-\epsilon \tau) + e^{-\pi{\tilde \sigma}/2} \Theta(\epsilon \tau) ) {\hat e}(\sigma,\epsilon,\tau) \otimes {\hat e}(s-\sigma+\tilde\sigma,\epsilon,t-\tau) + \right. \\
			& + \frac{e^{-\pi (s-\sigma)} (1-e^{-2\pi\sigma})}{1-e^{-2\pi s}} (e^{\pi{\tilde \sigma}/2} \Theta(\epsilon \tau) + e^{-\pi{\tilde \sigma}/2} \Theta(-\epsilon\tau) ) {\hat e}(\sigma,\epsilon,\tau) \otimes {\hat e}(s-\sigma+\tilde\sigma,-\epsilon,t-\tau) + \\
			&\left. + \frac{e^{-\pi\sigma}(1-e^{-2\pi(s-\sigma)})}{1-e^{-2\pi s}} (e^{\pi{\tilde \sigma}/2} \Theta(-\epsilon\tau) + e^{-\pi{\tilde \sigma}/2} \Theta(\epsilon \tau) ) {\hat e}(\sigma,-\epsilon,\tau) \otimes {\hat e}(s-\sigma+\tilde\sigma,\epsilon,t-\tau) \right\} .
	\end{aligned}\end{align}
\end{small}
The elements ${\hat e}(s,\epsilon,t)$ admit a presentation in terms of generators $\hat H, F$ satisfying equations \eqref{discrete_multiplication_dual}-\eqref{discrete_comultiplication_dual}
\begin{align}
	\hat e(s,\epsilon,t) &= |\hat{H}|^{is} \Theta(\epsilon \hat{H}) F^{ib^{-1}t}. \label{dualbasis}
\end{align}
As in the case of the initial Hopf algebra $\cA$, one can reproduce the discrete multiplication and co-multiplication relations by analytical continuation to the values $s=-im$, $t=-ibn$ of the powers, where $n,m\in\mathbb{Z}_{\geq0}$. \\

The algebras $\cA$ and $\cA^*$ are dual to each other in the sense of relation \eqref{duality} with respect to a duality pairing defined
\begin{align}
	(e(s,\epsilon,t),{\hat e}(s',\epsilon',t')) = \delta(s-s')\delta(t-t') \delta_{\epsilon,\epsilon'} .
\end{align}
Alternatively, one can see that the multiplication coefficients
\begin{align}\label{borelhalfmultiplication}
	\begin{aligned}
		e(s,\epsilon,t) e(s',\epsilon',t') &= \sum_{\epsilon''=\pm} \int\ud\sigma\ud\tau \, m_{s,\epsilon,t,s',\epsilon',t'}^{\sigma,\epsilon'',\tau} e(\sigma,\epsilon'',\tau), \\ 
		{\hat e}(s,\epsilon,t) {\hat e}(s',\epsilon',t') &= \sum_{\epsilon''=\pm} \int\ud\sigma\ud\tau \, {\hat m}^{s,\epsilon,t,s',\epsilon',t'}_{\sigma,\epsilon'',\tau} {\hat e}(\sigma,\epsilon'',\tau),
\end{aligned}\end{align}
and the co-multiplication coefficients
\begin{align}
	\label{borelhalfmultiplication2}
	\begin{aligned}
		\Delta(e(s,\epsilon,t)) &= \sum_{\epsilon',\epsilon''=\pm} \int\ud\sigma'\ud\sigma''\ud\tau'\ud\tau'' \, \mu_{s,\epsilon,t}^{\sigma',\epsilon',\tau',\sigma'',\epsilon'',\tau''} e(\sigma',\epsilon',\tau')\otimes e(\sigma'',\epsilon'',\tau''), \\ 
		{\hat \Delta}({\hat e}(s,\epsilon,t)) &= \sum_{\epsilon',\epsilon''=\pm} \int\ud\sigma'\ud\sigma''\ud\tau'\ud\tau'' \, {\hat \mu}^{s,\epsilon,t}_{\sigma',\epsilon',\tau',\sigma'',\epsilon'',\tau''} {\hat e}(\sigma',\epsilon',\tau')\otimes {\hat e}(\sigma'',\epsilon'',\tau''),
\end{aligned}\end{align}
defined as above do satisfy the following equalities
\begin{align}\begin{aligned}
		m_{\sigma',\epsilon',\tau',\sigma'',\epsilon'',\tau''}^{s,\epsilon,t} &= {\hat \mu}^{s,\epsilon,t}_{\sigma',\epsilon',\tau',\sigma'',\epsilon'',\tau''} , \\
		{\hat m}^{\sigma',\epsilon',\tau',\sigma'',\epsilon'',\tau''}_{s,\epsilon,t} &= \mu_{s,\epsilon,t}^{\sigma',\epsilon',\tau',\sigma'',\epsilon'',\tau''} .
\end{aligned}\end{align}
Given the above, one uses equation \eqref{heisenbergmultiplication} to define the commutation relations for the Heisenberg double of $\cA$. \\

We want to stress that, as before for $\cA$ and $\cA^*$ separately, one can reproduce the relations \eqref{commutationrelation} using the analytical continuation to the appropriate values of $\sigma, \tau$. When one analytically continues, the poles of gamma and Fadeev's quantum dilogarithm functions present in the integrand are pinching the contours of integration. The residues of those poles then contribute as the terms of the sums of the product of the discrete $H(\cA)$. \\

By applying the definition \eqref{canonicalelementdef} we obtain the following expression for the canonical element $S$ in terms of the basis elements of $H(\cA)$
\begin{equation}
	S= \sum_{\epsilon=\pm1} \int_{\mathbb{R}^2} \ud s\ud t \, e(s,\epsilon,t) \otimes {\hat e}(s,\epsilon,t) ,
\end{equation}
which can be written in terms of the generators $H, {\hat H}, E, F$ by using the explicit presentation of basis elements in equations \eqref{basis} and \eqref{dualbasis}, 
\begin{align}\label{ss}
	&S = e^{2\pi i H\otimes \hat H} g_{\ub}^{-1}(E \otimes F) ,
\end{align}
where $e_\ub(x) = g_\ub(e^{2 \pi b x})$. We see that quantum dilogarithm $(x,q)_\infty$, which was present in the discrete version of the algebra, has been replaced by its continuous analogue $g_\ub(x)$.

Moreover, the pentagon equation \eqref{pentagonrelation} for $S$ is manifestly satisfied as it can be reduced to the pentagon relation for Faddeev's quantum dilogarithm (equation \eqref{pentagonfaddeev} in appendix \ref{appendixA}).

\subsection*{Representations of the Heisenberg double of the Borel half of $U_q(sl(2))$}

In this part we will consider a representation theory of the continuous Heisenberg double $H(\mathcal{A})$ described above. This representations was first considered by Kashaev \cite{Ka1} in the context of applications to the \TM theory of Riemann surfaces. The Heisenberg double evaluated on those representations have a direct interpretation as the operators in the quantum \TM theory~\cite{Ka1,Ka2}. \\

Following \cite{Ka1}, we introduce the representations $\pi:H(\mathcal{A})\to Hom(L^2(\mathbb{R}))$ of the Heisenberg double $H(\mathcal{A})$ on $L^2(\mathbb{R})$ representation space using the following realisation of the generators $H,\hH,E,F$
\begin{align}\label{representation_sl2}
	\begin{aligned}
		&\pi(H) = \mathsf{p}, && \pi(E) = e^{2\pi b \mathsf{x}},\\
		&\pi(\hat H) = \mathsf{x}, && \pi(F) = e^{2\pi b(\mathsf{p}-\mathsf{x})},
	\end{aligned} 
\end{align}
where $\mathsf{p},\mathsf{x}$ are self-adjoint operators on $L^2(\mathbb{R})$ and $[\mathsf{p},\mathsf{x}]=\frac{1}{2\pi i}$. One can show that these generators  satisfy the commutation relations \eqref{commutationrelation}.

The canonical element $S$ \eqref{ss} evaluated on those representations can be written as follows
\begin{align}
	&(\pi\otimes\pi)(S) = e^{2\pi i \mathsf{p}_1 \mathsf{x}_2} e_{\ub}^{-1}(\mathsf{x}_1+\mathsf{p}_2-\mathsf{x}_2) .
\end{align}
This representation of the canonical element has been considered in the context of Teichm\"uller theory as a realisation of the flip operator \cite{Ka1,Ka2}.

\section{Countinuous Heisenberg double of the Borel half of $U_q(osp(1|2))$}\label{chapter4}

This section is devoted to the study of the continuous Heisenberg double of the Borel half of $U_q(osp(1|2))$, also known under a name of a quantum superplane. We construct the continuous Heisenberg double in a manner similar to the one described in section \ref{chapter3}, however the grading of elements is considered very carefully since $U_q(osp(1|2))$ is $\mathbb{Z}_2$-graded. Afterwards, we consider an infinite dimensional representations of the Heisenberg double on $L^2(\mathbb{R})\otimes\mathbb{C}^{1|1}$ with the focus on canonical element $S$. \\

The discrete Heisenberg double $H(\cA)$ of $\cA = \mathcal{B}(U_q(osp(1|2)))$ was already dissussed in the section \ref{chapter2}. We intend to find a continuous counterpart of that algebra, however this cannot be done in exactly the same fashion as the non-graded algebra like $\mathcal{B}(U_q(sl(2)))$. This is caused by the fact that one cannot take a complex power of elements which have an odd degree and produce a homogenous elements. In particular, we cannot simply take the imaginary powers of the odd elements $v^{(\pm)}$. In order to resolve this issue, we can consider a decomposition of those particular elements
\begin{align}\label{decomposition_osp12}\begin{aligned}
& v^{(+)} = E \kappa, ~~~~~~~~~~~~~
 &&v^{(-)} = F \hka ,
\end{aligned}\end{align}
into the even elements $E,F$ which satisfy non-trivial commutation relations with generators $H,\hH$
\begin{align}\begin{aligned}\label{decomposition_heisenberg_relations_osp12_eq1}
 & [H,E]=-ibE, &&  [\hat H,F]=+ib F.
\end{aligned}\end{align}
The odd elements $\ka,\hka$ commute trivially with the even ones
\begin{align}\begin{aligned}
 \phantom{yy} &[H,\ka] = [E,\ka] = 0, \\
 &[\hH,\hka] = [F,\hka] = 0, 
\end{aligned}\end{align}
and they satisfy the following identites
\begin{align}\label{decomposition_heisenberg_relations_osp12_eq3}
\begin{aligned}
 \kappa^2 = -1, && \hka^2 = -1.
 \end{aligned}
\end{align}
The decomposition \eqref{decomposition_osp12} informs one how one needs to modify the definition of the Hopf algebra $\cA$ in the continuous case. It allows to straightforwardly take the imaginary powers of the even part of the decomposition, while constraining the powers of the odd part to integers only. \\

Let us start 
with the Borel half $\cA$ of $U_q(osp(1|2))$. It is a Hopf algebra spanned by elements $\{e(s,\epsilon,t,n)\}_{s,t\in \mathbb{R}, \epsilon=\pm1,n=0,1}$, where the basis elements are given in terms of the generators by
\begin{align}\label{basis_osp12}
 e(s,\epsilon,t,n) &= f_{\epsilon,n}(s,t) ( \epsilon\pi H)^{is}_{pv} E^{ib^{-1}t} \kappa^n ,
\end{align}
where
\begin{align}\label{super-normalisation}\begin{aligned}
 f_{\epsilon,0}(s,t) &= \frac{1}{4\pi} \zeta_0 \Gamma(-is) e^{-\pi t Q/2} e^{\pi s/2} G_\NS^{-1}(Q+it) ,  \\
 f_{\epsilon,1}(s,t) &= \frac{i}{4\pi} \zeta_0 \Gamma(-is) e^{-\pi t Q/2} e^{\pi s/2} G_\R^{-1}(Q+it) ,
\end{aligned}\end{align}
where the special functions $G_\R, G_\NS$ are related to the supersymmetric analogues of Faddeev's quantum dilogarithm functions
\begin{align}\label{super_quantum_dilogs}\begin{aligned}
 e^{-1}_\NS(r) + e^{-1}_\R(r) &= \zeta_0 \int\ud t e^{ \pi itr} \frac{e^{-\frac{1}{2}\pi tQ}}{G_\NS(Q+it)}  , \\
 e^{-1}_\NS(r) - e^{-1}_\R(r) &= \zeta_0 \int\ud t e^{ \pi itr} \frac{e^{-\frac{1}{2}\pi tQ}}{G_\R(Q+it)}  ,
\end{aligned}\end{align}
that are described in more details in the appendix \ref{appendixA}, with $\zeta_0 = \exp(-i\pi (b+b^{-1})^2/8)$. 

The basis elements have the following product relations
 \begin{align}\begin{aligned}\label{continuous_multiplication_osp12_1}
 & e(s,\epsilon,t,n) e(s',\epsilon,t',n')= \int \frac{\ud \sigma}{2\pi} {{-is'}\choose{-i\sigma}}_{\Gamma} |\pi t|^{i\sigma} \frac{f_{\epsilon,n}(s,t)f_{\epsilon,n'}(s',t')}{f_{\epsilon,n+n'}(s+s'-\sigma,t+t')} \times \\
 &\times [\Theta(-\epsilon t) +\Theta(\epsilon t)e^{-\pi \sigma }] e(s+s'-\sigma,\epsilon,t+t',n+n') ,\\
 & e(s,\epsilon,t,n) e(s',-\epsilon,t',n') = \int \frac{\ud \sigma}{2\pi} {{-is'}\choose{-i\sigma}}_{\Gamma} |\pi t|^{i\sigma} [\Theta(\epsilon t) +\Theta(-\epsilon t)e^{-\pi \sigma }] \times \\
&\times \left[ \frac{e^{-\pi(s'-\sigma)}(1-e^{-2\pi s})}{1-e^{-2 \pi(s'+s-\sigma)}} \frac{f_{\epsilon,n}(s,t)f_{-\epsilon,n'}(s',t')}{f_{\epsilon,n+n'}(s+s'-\sigma,t+t')} e(s+s'-\sigma,\epsilon,t+t',n+n') + \right. \\ 
&\left. + \frac{e^{-\pi s}(1-e^{-2\pi (s'-\sigma)})}{1-e^{-2 \pi(s'+s-\sigma)}} \frac{f_{\epsilon,n}(s,t)f_{-\epsilon,n'}(s',t')}{f_{-\epsilon,n+n'}(s+s'-\sigma,t+t')} e(s+s'-\sigma,-\epsilon,t+t',n+n') \right],
\end{aligned}\end{align}
while the co-product is given by 
\begin{align}\begin{footnotesize}\begin{aligned}
& \Delta(e(s,\epsilon,t,1)) = \frac{1}{2} \zeta_0\int\frac{\ud\sigma \ud{\tilde \sigma}}{(2\pi)^2}\ud\tau \Gamma(-i\tilde\sigma) { -is \choose -i\sigma}_\Gamma e^{\pi{\tilde\sigma}/2} |\tau|^{i\tilde\sigma} \left[ \frac{G_\R(Q+it)}{G_\NS(Q+i\tau) G_\R(-i\tau+Q+it)} \times \right. \\
&\times e(\sigma,\epsilon,\tau,0) \otimes \Bigg\{ \left( \Theta(\epsilon\tau) + \frac{e^{-\pi\tilde\sigma}(1-e^{-2\pi(s-\sigma)})}{1-e^{-2\pi(s-\sigma+\tilde\sigma)}} \Theta(-\epsilon \tau) \right) \frac{f_{\epsilon,1}(s,t)}{f_{\epsilon,0}(\sigma,\tau) f_{\epsilon,1}(s-\sigma+\tilde\sigma,t-\tau)} \times \\
&\left. \times e(s-\sigma+\tilde\sigma,\epsilon,t-\tau,1) +\frac{e^{-\pi(s-\sigma)}(1-e^{-2\pi\tilde\sigma})}{1-e^{-2\pi(s-\sigma+\tilde\sigma)}} \Theta(-\epsilon \tau) \frac{f_{\epsilon,1}(s,t)}{f_{\epsilon,0}(\sigma,\tau) f_{-\epsilon,1}(s-\sigma+\tilde\sigma,t-\tau)} \times \right. \\
& \times e(s-\sigma+\tilde\sigma,-\epsilon,t-\tau,1) \Bigg\}+ \frac{G_\R(Q+it)}{G_\R(Q+i\tau)G_\NS(-i\tau+Q+it)} e(\sigma,\epsilon,\tau,1) \otimes \\
&\otimes \left\{ \left( \Theta(\epsilon\tau) + \frac{e^{-\pi\tilde\sigma}(1-e^{-2\pi(s-\sigma)})}{1-e^{-2\pi(s-\sigma+\tilde\sigma)}} \Theta(-\epsilon\tau) \right) \frac{f_{\epsilon,1}(s,t)}{f_{\epsilon,1}(\sigma,\tau) f_{\epsilon,0}(s-\sigma+\tilde\sigma,t-\tau)} e(s-\sigma+\tilde\sigma,\epsilon,t-\tau,0) + \right. \\
&\left.\left. +\frac{e^{-\pi(s-\sigma)}(1-e^{-2\pi\tilde\sigma})}{1-e^{-2\pi(s-\sigma+\tilde\sigma)}} \Theta(-\epsilon\tau) \frac{f_{\epsilon,1}(s,t)}{f_{\epsilon,1}(\sigma,\tau) f_{-\epsilon,0}(s-\sigma+\tilde\sigma,t-\tau)} e(s-\sigma+\tilde\sigma,-\epsilon,t-\tau,0) \right\} \right], 
\end{aligned}\end{footnotesize}\end{align}
\begin{align}\begin{footnotesize}\begin{aligned}\label{continuous_multiplication_osp12_3}
& \Delta(e(s,\epsilon,t,0)) = \frac{1}{2} \zeta_0\int\frac{\ud\sigma \ud{\tilde \sigma}}{(2\pi)^2}\ud\tau \Gamma(-i\tilde\sigma) { -is \choose -i\sigma}_\Gamma e^{\pi{\tilde\sigma}/2} |\tau|^{i\tilde\sigma} \left[ \frac{G_\NS(Q+it)}{G_\NS(Q+i\tau) G_\NS(-i\tau+Q+it)} \times \right. \\
&\times e(\sigma,\epsilon,\tau,0) \otimes \Bigg\{ \left( \Theta(\epsilon\tau) + \frac{e^{-\pi\tilde\sigma}(1-e^{-2\pi(s-\sigma)})}{1-e^{-2\pi(s-\sigma+\tilde\sigma)}} \Theta(-\epsilon \tau) \right) \frac{f_{\epsilon,0}(s,t)}{f_{\epsilon,0}(\sigma,\tau) f_{\epsilon,0}(s-\sigma+\tilde\sigma,t-\tau)} \times \\
&\left. \times e(s-\sigma+\tilde\sigma,\epsilon,t-\tau,0) +\frac{e^{-\pi(s-\sigma)}(1-e^{-2\pi\tilde\sigma})}{1-e^{-2\pi(s-\sigma+\tilde\sigma)}} \Theta(-\epsilon \tau) \frac{f_{\epsilon,0}(s,t)}{f_{\epsilon,0}(\sigma,\tau) f_{-\epsilon,0}(s-\sigma+\tilde\sigma,t-\tau)} \times \right. \\
& \times e(s-\sigma+\tilde\sigma,-\epsilon,t-\tau,0) \Bigg\}+ \frac{G_\NS(Q+it)}{G_\R(Q+i\tau)G_\R(-i\tau+Q+it)} e(\sigma,\epsilon,\tau,1) \otimes \\
&\otimes \left\{ \left( \Theta(\epsilon\tau) + \frac{e^{-\pi\tilde\sigma}(1-e^{-2\pi(s-\sigma)})}{1-e^{-2\pi(s-\sigma+\tilde\sigma)}} \Theta(-\epsilon\tau) \right) \frac{f_{\epsilon,0}(s,t)}{f_{\epsilon,1}(\sigma,\tau) f_{\epsilon,1}(s-\sigma+\tilde\sigma,t-\tau)} e(s-\sigma+\tilde\sigma,\epsilon,t-\tau,1) + \right. \\
&\left.\left. +\frac{e^{-\pi(s-\sigma)}(1-e^{-2\pi\tilde\sigma})}{1-e^{-2\pi(s-\sigma+\tilde\sigma)}} \Theta(-\epsilon\tau) \frac{f_{\epsilon,0}(s,t)}{f_{\epsilon,1}(\sigma,\tau) f_{-\epsilon,1}(s-\sigma+\tilde\sigma,t-\tau)} e(s-\sigma+\tilde\sigma,-\epsilon,t-\tau,1) \right\} \right].
\end{aligned}\end{footnotesize}\end{align}

By analytically continuing the values of $s, t$ in the equations \eqref{continuous_multiplication_osp12_1}-\eqref{continuous_multiplication_osp12_3} one can recover the commutation \eqref{discrete_osp_multiplication_basis} and co-product relations \eqref{discrete_osp_comultiplication_basis} for the discrete algebra elements generated by $H, \vp$ considered in section \ref{chapter2}. The values corresponding to the discrete algebra basis element are $s=-im$, $t=-ib n$, $n=1$ for $m\in\mathbb{Z}_{\geq0},n\in2\mathbb{Z}_{\geq0}+1$ and $s=-im$, $t=-ibn$, $n=0$ for $m\in\mathbb{Z}_{\geq0},n\in2\mathbb{Z}_{\geq0}$. \\

After describing the Hopf algebra $\cA$, we consider the dual Hopf algebra $\cA^*$. This Hopf algebra is spanned by the elements $\{\hat e(s,\epsilon,t,n)\}_{s,t\in \mathbb{R}, \epsilon=\pm1, n=0,1}$. The dual basis can be expressed in terms of the generators satisfying the relations \eqref{decomposition_heisenberg_relations_osp12_eq1}-\eqref{decomposition_heisenberg_relations_osp12_eq3}
\begin{align}\label{dualbasis_osp12}
 \hat e(s,\epsilon,t,n) &= |\hat{H}|^{is} \Theta(\epsilon \hat{H}) F^{ib^{-1}t} \hka^n .
 \end{align}
The multiplication relations for those elements are as follows
 \begin{small}\begin{align}\begin{aligned}
 & {\hat e}(s,\epsilon,t,n){\hat e}(s',\epsilon,t',n') = \int \frac{\ud \sigma}{2\pi} {{-is'}\choose{-i\sigma}}_{\Gamma} |t|^{i\sigma} \left[\Theta(-\epsilon t) \frac{e^{-\pi \sigma}(1-e^{-2\pi (s'-\sigma)})}{1-e^{-2\pi s'}}+\Theta({\epsilon}t) \right] \times \\
 &\times {\hat e}(s+s'-\sigma,\epsilon,t+t',n+n'),\\
 & {\hat e}(s,\epsilon,t,n){\hat e}(s',-\epsilon,t',n') = \int \frac{\ud \sigma}{2\pi} {{-is'}\choose{-i\sigma}}_{\Gamma} |t|^{i\sigma} \Theta(-\epsilon t) \frac{e^{-\pi (s'-\sigma)}(1-e^{-2\pi \sigma})} {1-e^{-2\pi s'}} \times \\ 
 &\times {\hat e}(s+s'-\sigma,\epsilon,t+t',n+n') ,
 \end{aligned}\end{align}                         \end{small}
while the co-multiplication has the following form 
\begin{align}\begin{footnotesize}\begin{aligned}
&\Delta({\hat e}(s,\epsilon,t,0)) = \frac{1}{2} \zeta_0 \int \frac{\ud\sigma}{2\pi} \frac{\ud {\tilde \sigma}}{2\pi} \ud\tau \Gamma(-i{\tilde \sigma}) {{-is}\choose{-i\sigma}}_{\Gamma}  |\pi\tau|^{i{\tilde \sigma}} \times \\
&\times \left[ \frac{G_\NS(Q+it)}{G_\NS(Q+i\tau) G_\NS(-i\tau+Q+it)} \left\{ (e^{\pi{\tilde \sigma}/2} \Theta(-\epsilon\tau) + e^{-\pi{\tilde \sigma}/2} \Theta(\epsilon\tau) ) {\hat e}(\sigma,\epsilon,\tau,0) \otimes {\hat e}(s-\sigma+\tilde\sigma,\epsilon,t-\tau,0) + \right.\right. \\
& + \frac{e^{-\pi (s-\sigma)} (1-e^{-2\pi\sigma})}{1-e^{-2\pi s}} (e^{\pi{\tilde \sigma}/2} \Theta(\epsilon\tau) + e^{-\pi{\tilde \sigma}/2} \Theta(-\epsilon\tau) ) {\hat e}(\sigma,\epsilon,\tau,0) \otimes {\hat e}(s-\sigma+\tilde\sigma,-\epsilon,t-\tau,0) + \\
&\left. + \frac{e^{-\pi\sigma}(1-e^{-2\pi(s-\sigma)})}{1-e^{-2\pi s}} (e^{\pi{\tilde \sigma}/2} \Theta(-\epsilon\tau) + e^{-\pi{\tilde \sigma}/2} \Theta(\epsilon\tau) ) {\hat e}(\sigma,-\epsilon,\tau,0) \otimes {\hat e}(s-\sigma+\tilde\sigma,\epsilon,t-\tau,0) \right\} + \\
& + \frac{G_\NS(Q+it)}{G_\R(Q+i\tau) G_\R(-i\tau+Q+it)} \left\{  (e^{\pi{\tilde \sigma}/2} \Theta(-\epsilon\tau) + e^{-\pi{\tilde \sigma}/2} \Theta(\epsilon\tau) ) {\hat e}(\sigma,\epsilon,\tau,1) \otimes {\hat e}(s-\sigma+\tilde\sigma,\epsilon,t-\tau,1) + \right. \\
& + \frac{e^{-\pi (s-\sigma)} (1-e^{-2\pi\sigma})}{1-e^{-2\pi s}} (e^{\pi{\tilde \sigma}/2} \Theta(\epsilon\tau) + e^{-\pi{\tilde \sigma}/2} \Theta(-\epsilon\tau) ) {\hat e}(\sigma,\epsilon,\tau,1) \otimes {\hat e}(s-\sigma+\tilde\sigma,-\epsilon,t-\tau,1) + \\
&\left.\left. + \frac{e^{-\pi\sigma}(1-e^{-2\pi(s-\sigma)})}{1-e^{-2\pi s}} (e^{\pi{\tilde \sigma}/2} \Theta(-\epsilon\tau) + e^{-\pi{\tilde \sigma}/2} \Theta(\epsilon\tau) ) {\hat e}(\sigma,-\epsilon,\tau,1) \otimes {\hat e}(s-\sigma+\tilde\sigma,\epsilon,t-\tau,1) \right\} \right] , 
\end{aligned}\end{footnotesize}\end{align}
\begin{align}\begin{footnotesize}\begin{aligned}
&\Delta({\hat e}(s,\epsilon,t,1)) = \frac{1}{2} \zeta_0 \int \frac{\ud\sigma}{2\pi} \frac{\ud {\tilde \sigma}}{2\pi} \ud\tau \Gamma(-i{\tilde \sigma}) {{-is}\choose{-i\sigma}}_{\Gamma}  |\pi\tau|^{i{\tilde \sigma}} \times \\
&\times \left[ \frac{G_\R(Q+it)}{G_\R(Q+i\tau) G_\NS(-i\tau+Q+it)} \left\{ (e^{\pi{\tilde \sigma}/2} \Theta(-\epsilon\tau) + e^{-\pi{\tilde \sigma}/2} \Theta(\epsilon\tau) ) {\hat e}(\sigma,\epsilon,\tau,1) \otimes {\hat e}(s-\sigma+\tilde\sigma,\epsilon,t-\tau,0) + \right.\right. \\
& + \frac{e^{-\pi (s-\sigma)} (1-e^{-2\pi\sigma})}{1-e^{-2\pi s}} (e^{\pi{\tilde \sigma}/2} \Theta(\epsilon\tau) + e^{-\pi{\tilde \sigma}/2} \Theta(-\epsilon\tau) ) {\hat e}(\sigma,\epsilon,\tau,1) \otimes {\hat e}(s-\sigma+\tilde\sigma,-\epsilon,t-\tau,0) + \\
&\left. + \frac{e^{-\pi\sigma}(1-e^{-2\pi(s-\sigma)})}{1-e^{-2\pi s}} (e^{\pi{\tilde \sigma}/2} \Theta(-\epsilon\tau) + e^{-\pi{\tilde \sigma}/2} \Theta(\epsilon\tau) ) {\hat e}(\sigma,-\epsilon,\tau,1) \otimes {\hat e}(s-\sigma+\tilde\sigma,\epsilon,t-\tau,0) \right\} + \\
& + \frac{G_\R(Q+it)}{G_\NS(Q+i\tau) G_\R(-i\tau+Q+it)} \left\{  (e^{\pi{\tilde \sigma}/2} \Theta(-\epsilon\tau) + e^{-\pi{\tilde \sigma}/2} \Theta(\epsilon\tau) ) {\hat e}(\sigma,\epsilon,\tau,0) \otimes {\hat e}(s-\sigma+\tilde\sigma,\epsilon,t-\tau,1) + \right. \\
& + \frac{e^{-\pi (s-\sigma)} (1-e^{-2\pi\sigma})}{1-e^{-2\pi s}} (e^{\pi{\tilde \sigma}/2} \Theta(\epsilon\tau) + e^{-\pi{\tilde \sigma}/2} \Theta(-\epsilon\tau) ) {\hat e}(\sigma,\epsilon,\tau,0) \otimes {\hat e}(s-\sigma+\tilde\sigma,-\epsilon,t-\tau,1) + \\
&\left.\left. + \frac{e^{-\pi\sigma}(1-e^{-2\pi(s-\sigma)})}{1-e^{-2\pi s}} (e^{\pi{\tilde \sigma}/2} \Theta(-\epsilon\tau) + e^{-\pi{\tilde \sigma}/2} \Theta(\epsilon\tau) ) {\hat e}(\sigma,-\epsilon,\tau,0) \otimes {\hat e}(s-\sigma+\tilde\sigma,\epsilon,t-\tau,1) \right\} \right] .
\end{aligned}\end{footnotesize}\end{align} 
As in the case of $\cA$, the multiplication and co-multiplication relations above reduce to the product and co-product \eqref{discrete_osp_multiplication_basis_dual}-\eqref{discrete_osp_comultiplication_basis_dual} of the discrete dual Hopf algebra from section \ref{chapter2} by the means of appropriate analytic continuation. \\

One can see that the multiplications and co-multiplications of $\cA$ and $\cA^*$ are dual to each other in the sense of \eqref{duality} with respect to a duality bracked defined
\begin{align}
 (e(s,\epsilon,t,n),{\hat e}(s',\epsilon',t',n')) = \delta(s-s')\delta(t-t') \delta_{\epsilon,\epsilon'} \delta_{n,n'} .
\end{align}
Alternatively, one can see that with the multiplication and co-multiplication coefficients defined in the following way
\begin{small}\begin{align}\label{borelhalfmultiplication_osp12}
	\begin{aligned}
 e(s,\epsilon,t,n) e(s',\epsilon',t',n') &= \sum_{\epsilon''=\pm}\sum_{n''=0}^1 \int\ud\sigma\ud\tau \, m_{s,\epsilon,t,n,s',\epsilon',t',n'}^{\sigma,\epsilon'',\tau,n''} e(\sigma,\epsilon'',\tau,n''), \\ 
 {\hat e}(s,\epsilon,t,n) {\hat e}(s',\epsilon',t',n') &= \sum_{\epsilon''=\pm}\sum_{n''=0}^1 \int\ud\sigma\ud\tau \, {\hat m}^{s,\epsilon,t,n,s',\epsilon',t',n'}_{\sigma,\epsilon'',\tau,n''} {\hat e}(\sigma,\epsilon'',\tau,n''),
\end{aligned}\end{align}
\end{small}and
\begin{small}\begin{align}
	\label{borelhalfmultiplication2_osp12}
	\begin{aligned}
 \Delta(e(s,\epsilon,t,n)) &= \sum_{\substack{\epsilon',\epsilon''=\pm,\\n',n''=0,1}} \int\ud\sigma'\ud\sigma''\ud\tau'\ud\tau'' \, \mu_{s,\epsilon,t,n}^{\sigma',\epsilon',\tau',n',\sigma'',\epsilon'',\tau'',n''} e(\sigma',\epsilon',\tau',n')\otimes e(\sigma'',\epsilon'',\tau'',n''), \\ 
 {\hat \Delta}({\hat e}(s,\epsilon,t,n)) &= \sum_{\substack{\epsilon',\epsilon''=\pm,\\n',n''=0,1}} \int\ud\sigma'\ud\sigma''\ud\tau'\ud\tau'' \, {\hat \mu}^{s,\epsilon,t,n}_{\sigma',\epsilon',\tau',n',\sigma'',\epsilon'',\tau'',n''} {\hat e}(\sigma',\epsilon',\tau',n')\otimes {\hat e}(\sigma'',\epsilon'',\tau'',n''),
\end{aligned}\end{align}
\end{small}satisfy the equality
\begin{align}\begin{aligned}
  m_{\sigma',\epsilon',\tau',n',\sigma'',\epsilon'',\tau'',n''}^{s,\epsilon,t,n} &= (-1)^{n' n''} {\hat \mu}^{s,\epsilon,t,n}_{\sigma',\epsilon',\tau',n',\sigma'',\epsilon'',\tau'',n''} , \\
  {\hat m}^{\sigma',\epsilon',\tau',n',\sigma'',\epsilon'',\tau'',n''}_{s,\epsilon,t,n} &= (-1)^{n' n''} \mu_{s,\epsilon,t,n}^{\sigma',\epsilon',\tau',n',\sigma'',\epsilon'',\tau'',n''} .
\end{aligned}\end{align}
Then, one uses \eqref{heisenbergmultiplication} to find the exchange relations for the Heisenberg double $H(\cA)$. By analytic continuation, one can produce the (anti-)commutation relations for the generators $H,\hat H, E, F, \kappa, \hat\kappa$ which are as follows: the even generators satisfy 
\begin{align}\begin{aligned}\label{heisenberg_relations_osp12_eq1}
 & [E, F] = e^{\pi bH}(q^{\frac{1}{2}}-q^{-\frac{1}{2}}) ,\\
 & [H,E]=-ibE, && [H, F] = +ib F, \\
 & [\hat H, E] = 0, && [\hat H,F]=+ib F,
\end{aligned}\end{align}
while $\ka,\hka$ commute trivially with all even generators
\begin{align}\begin{aligned}
 \phantom{yy} &[H,\ka] = [\hH,\ka] = [E,\ka] = [F,\ka] = 0, \\
 &[H,\hka] = [\hH,\hka] = [E,\hka] = [F,\hka] = 0, 
\end{aligned}\end{align}
and they satisfy the following identites between each other
\begin{align}\label{heisenberg_relations_osp12_eq3}
\begin{aligned}
 &\kappa^2 = \hka^2 = -1, ~~~~~~~~~~~~~&&\ka = \hka .
 \end{aligned}
\end{align}
It is worthwile to note that, taking into account the decomposition \eqref{decomposition_osp12} from which we started this section, one can recover the (anti-)commutation relations of the discrete Heisenberg double \eqref{discrete_osp_heisenberg_commutation_relations} using the exchange relations that come from \eqref{heisenbergmultiplication}. \\

Let us now take a look at the canonical element $S$ of this Heisenberg double. Using the definition \eqref{canonicalelementdef} we obtain the relation for $S$ in terms of the basis elements \eqref{basis_osp12} and \eqref{dualbasis_osp12}
\begin{equation}
 S= \sum_{\epsilon=\pm1} \sum_{n=0}^1 \int \ud s\ud t \, (-1)^{n} e(s,\epsilon,t,n) \otimes {\hat e}(s,\epsilon,t,n) ,
\end{equation}
which written in terms of the generators has the form
\begin{equation}\begin{aligned}\label{canonicalelement_osp}
 S &= \frac{1}{2} e^{i \pi H \otimes {\hat H}} \left[ \left( g_\R^{-1}(E\otimes F) + g_\NS^{-1}(E\otimes F) \right) 1\otimes 1 +\right. \\
 & \left. \qquad +i \left( g_\R^{-1}(E\otimes F) - g_\NS^{-1}(E\otimes F) \right) \kappa \otimes {\hat \kappa} \right] ,
\end{aligned}\end{equation}
where we used equation \eqref{super_quantum_dilogs} and the relations
\begin{align}
 & g_\R(x) = e_\R(e^{\pi b x}), && g_\NS(x) = e_\NS(e^{\pi b x}) .
\end{align}
The fact that the canonical element $S$, equation \eqref{canonicalelement_osp}, satisfies the graded pentagon equation follows directly from the supersymmetric pentagon identities \eqref{superpentagon} and has been checked explicitly in \cite{APT}.

\subsection*{Representations of the Heisenberg double of the Borel half of $U_q(osp(1|2))$}

In this section we want to introduce the representation $\pi:H(\cA)\to Hom(L^2(\mathbb{R})\otimes\mathbb{C}^{1|1})$ of the Heisenberg double $H(\cA)$ of the quantum superplane $\cA$ that is a supersymmetric analogue of the representation \eqref{representation_sl2}. The generators are represented as the following operators
\begin{align}\label{representation_osp12}
\begin{aligned}
 & \pi(H) = \mathsf{p} \left( \begin{array}{cc} 1 & 0\\0& 1 \end{array} \right), && \pi(E) = e^{\pi b \mathsf{x}} \left( \begin{array}{cc} 1 & 0\\0& 1 \end{array} \right), && \pi(\ka) = i \left( \begin{array}{cc} 0 & 1\\1& 0 \end{array} \right),\\ 
 &\pi(\hat H) = \mathsf{x} \left( \begin{array}{cc} 1 & 0\\0& 1 \end{array} \right), && \pi(F) = e^{\pi b(\mathsf{p}-\mathsf{x})} \left( \begin{array}{cc} 1 & 0\\0& 1 \end{array} \right), && \pi(\hka) = i \left( \begin{array}{cc} 0 & 1\\1& 0 \end{array} \right),
 \end{aligned}
\end{align}
where $[\mathsf{p},\mathsf{x}]=\frac{1}{\pi i}$ are operators on $L^2(\mathbb{R})$. The canonical element $S$ \eqref{canonicalelement_osp} evaluated on the representation \eqref{representation_osp12} has the form
\begin{equation}
 \begin{split}
  (\pi\otimes\pi)(S) &= \frac{1}{2} e^{\pi i \mathsf{p}_1 \mathsf{x}_2} \Bigg\{ \Big[e_\R^{-1}(\mathsf{x}_1+\mathsf{p}_2-\mathsf{x}_2)+e_\NS^{-1}(\mathsf{x}_1+\mathsf{p}_2-\mathsf{x}_2)\Big] \left( \begin{array}{cc} 1 & 0\\0& 1 \end{array} \right) \otimes \left( \begin{array}{cc} 1 & 0\\0& 1 \end{array} \right) + \\
 &\qquad - i \Big[e_\R^{-1}(\mathsf{x}_1+\mathsf{p}_2-\mathsf{x}_2)-e_\NS^{-1}(\mathsf{x}_1+\mathsf{p}_2-\mathsf{x}_2)\Big] \left( \begin{array}{cc} 0 & 1\\1& 0 \end{array} \right) \otimes \left( \begin{array}{cc} 0 & 1\\1& 0 \end{array} \right) \Bigg\} .
 \end{split}
\end{equation}

This representation of the canonical element has been considered in the context of super Teichm\"uller theory as a realisation of the supersymmetric flip operator \cite{APT}.

\section{Drinfeld double}\label{chapter5}

In this section, we will present the definition of the $\mathbb{Z}_2$-graded Drinfeld double $D(\cA)$ of a Hopf algebra $\cA$, given in terms of the basis elements, and remind ourselves some facts about the universal element $R$ satisfying the Yang-Baxter equation. Then, we will describe an algebra morphism between $H(\cA)\otimes H(\cA^*)$ and $D(\cA)$, which constitutes a $\mathbb{Z}_2$-graded generalisation of a morphism described in \cite{Ka3}. Furthermore, we state the relation between the universal elements of the Heisenberg doubles and the universal $R$-matrix of the Drinfeld double. We present the $U_q(sl(2))$ and $U_q(osp(1|2))$ algebras as instructive examples. \\

Let us consider again a Hopf algebra $(\mathcal{A},m,\eta,\Delta,\epsilon,\gamma)$, subjected to the axioms \eqref{hopf_alg_axioms_1}-\eqref{hopf_alg_axioms_3}. With the choice of a basis $\{e_\alpha\}_{\alpha\in I}$ which algebraically spans $\cA$ we can describe the multiplication and co-multiplication 
\begin{align}
& e_\alpha e_\beta = \sum_{\gamma\in I} m^\gamma_{\alpha\beta} e_\gamma,
&& \Delta(e_\alpha) = \sum_{\beta, \gamma\in I} \mu^{\beta\gamma}_\alpha e_\beta\otimes e_\gamma ,
\end{align}
and an antipode 
\begin{align}
 \gamma(e_\alpha) = \sum_{\beta\in I} \gamma_\alpha^\beta e_\beta .
\end{align}
In addition to the algebra $\cA$ we can consider the algebra $\cA^*$, that is a Hopf algebra dual to $\cA$. In terms of a basis $\{e^\alpha\}_{\alpha\in I}$ of $\cA^*$ the multiplication and co-multiplication relations of this Hopf algebra are as follows
\begin{align}
&e^\alpha e^\beta = \sum_{\gamma\in I} (-1)^{|\alpha||\beta|} \mu_\gamma^{\alpha\beta} e^\gamma,
&&{\hat \Delta}^\text{op}(e^\alpha) = \sum_{\beta, \gamma\in I} (-1)^{|\beta||\gamma|} m_{\beta \gamma}^\alpha e^\beta\otimes e^\gamma ,
\end{align}
with an antipode 
\begin{align}
 {\hat \gamma}(e^\alpha) = \sum_{\beta\in I} \gamma^\alpha_\beta e^\beta .
\end{align}
They are dual to each other with respect to a duality bracket $(,)$ satisfying the relations \eqref{duality_bracket_eq1}-\eqref{duality_bracket_eq3}, and given explicitly on the basis by \eqref{duality}.

Given those two Hopf algebras, it is possible to define a quasi-triangular Hopf algebra as a double cross product of Hopf algebras \cite{Majidbook}. Explicitly, we can define a Hopf algebra $D(\cA)=\{ x\otimes f | x\in\cA, f\in\cA^*\}$ equipped with the product 
\begin{equation}\label{drinfeld_definition_eq1}\begin{aligned}
 (x\otimes (-1)^{|f|} f)(y\otimes (-1)^{|g|} g) &= \sum_{(y), (f)} (-1)^{|y_{(1)}| |f| + |y_{(2)}| (|f_{(1)}|+|f_{(2)}|) + |y_{(3)}| |f_{(1)}| + |f_{(2)}| |f_{(3)}| } \times \\
 &\qquad\times (-1)^{|f_{(1)}| (|f_{(2)}|+|f_{(3)}|) + |x| (|y_{(1)}| + |f_{(3)}|) + |g| (|f_{(1)}| + |x_{(3)}|)} \times \\
 &\qquad\times ( y_{(1)}, S^{-1}(f_{(3)}) ) ( y_{(3)}, f_{(1)} ) x y_{(2)} \otimes (-1)^{|f_{(2)}| + |g|} f_{(2)} g ,
\end{aligned}\end{equation}
the co-product
\begin{align}
 \Delta(x\otimes f) = \sum_{(x), (f)} (-1)^{|f_{(1)}||f_{(2)}|+|x_{(2)}||f_{(2)}|} x_{(1)} \otimes f_{(2)} \otimes x_{(2)} \otimes f_{(1)} , 
\end{align}
and the antipode
\begin{align}\label{drinfeld_definition_eq3}
 \gamma(x\otimes f) = (-1)^{|x||f|} (1 \otimes {\hat \gamma}^\text{cop}(f))(\gamma(x)\otimes 1),
\end{align}
where ${\hat \gamma}^\text{cop} = {\hat \gamma}^{-1}$ is the antipode of $(\mathcal{A}^*)^\text{cop}$, i.e. the antipode of the Hopf algebra dual to $\cA$ equipped with the opposite coproduct. We also use the following notation for the coproduct
\begin{align*}
 (\Delta\otimes id)\Delta(x) = \sum_{(x)} x_{(1)} \otimes x_{(2)} \otimes x_{(3)}, && ({\hat \Delta}\otimes id){\hat \Delta}(f) = \sum_{(f)} f_{(1)} \otimes f_{(2)} \otimes f_{(3)} .
\end{align*}

\begin{defn}
A {\normalfont{Drinfeld double}} of a Hopf algebra $\cA$ is a quasi-triangular Hopf algebra $D(\cA)$ with the multiplication, co-multiplication and antipode given by the equations \eqref{drinfeld_definition_eq1}-\eqref{drinfeld_definition_eq3}.
\end{defn}

The Hopf algebra $D(\cA)$ has $\mathcal{A}$ and $(\mathcal{A}^*)^\text{cop}$ as subalgebras through canonical embeddings $(\mathcal{A}^*)^\text{cop} \ni f\mapsto 1\otimes f\in D(\mathcal{A})$ and $\mathcal{A}\ni x \mapsto x\otimes 1 \in D(\mathcal{A})$. \\

In terms of a basis, the Drinfeld double $D(\cA)$ is spanned by a collection of elements $\{e_\alpha \otimes e^\beta \}_{\alpha,\beta\in I}$ which satisfy the following multiplication relations
\begin{equation}\label{susyEE}
\begin{split}
& (e_\alpha \otimes1) (1\otimes e^\beta) = e_\alpha \otimes e^\beta ,\\
& (e_\alpha \otimes1) (e_\beta \otimes 1) = m^\gamma_{\alpha\beta} (e_\gamma\otimes 1),\\
& (1\otimes e^\alpha) (1\otimes e^\beta) = (-1)^{|\alpha||\beta|}\mu_\gamma^{\alpha\beta} (1\otimes e^\gamma), \\
& (1\otimes e^{\alpha} ) (e_{\beta} \otimes 1)= (-1)^{|\mu|(|\sigma|+|\delta|)}  m_{\nu\gamma}^{\mu}m_{\mu\delta}^{\alpha} \mu_{\beta}^{\epsilon \rho}\mu_{\rho}^{\sigma \nu} ({\gamma}^{-1})^{\delta}_{\epsilon} e_{\sigma} \otimes e^{\gamma} ,
\end{split}
\end{equation}
and co-multiplication relations
\begin{align}
  \Delta(e_\alpha\otimes1) = \sum_{\beta, \gamma\in I} \mu^{\beta\gamma}_\alpha (e_\beta\otimes1)\otimes (e_\gamma\otimes1), &&   \Delta(1\otimes e^\alpha) = \sum_{\beta, \gamma\in I} m_{\gamma\beta}^\alpha (1\otimes e^\beta)\otimes (1\otimes e^\gamma) ,
\end{align}
and is equipped with the antipode
\begin{align}
  \gamma(e_\alpha\otimes 1) = \sum_{\beta\in I} \gamma_\alpha^\beta e_\beta\otimes 1, && \gamma(1\otimes e^\alpha) = \sum_{\beta\in I} {(\gamma^\text{cop})}_\alpha^\beta 1\otimes e^\beta ,
\end{align}
where $(\hat\gamma^\text{cop})^\alpha_\beta=(\gamma^{-1})^\alpha_\beta$. Equivalently, instead of the 4th exchange relation in \eqref{susyEE} one can use the crossing relation
\begin{align}
&\sum_{\gamma,\rho,\sigma\in I} (-1)^{|\beta||\sigma|} \mu^{\sigma\gamma}_{\alpha} m^\beta_{\gamma\rho} e_\sigma \otimes e^\rho = \sum_{\gamma,\rho,\sigma\in I} (-1)^{|\rho||\gamma|} m^\beta_{\rho\gamma} \mu^{\gamma\sigma}_\alpha (1\otimes e^\rho) \otimes ( e_\sigma \otimes 1) ,\label{susytt}
\end{align}
which indeed defines the same algebra. \\

\begin{defn}
In the case of Drinfeld double one has a {\normalfont{canonical element}} $R \in D(\mathcal{A})\otimes D(\mathcal{A})$ called the {\normalfont{universal $R$-matrix}}
\begin{equation}\label{rmatrix_def}
 R = \sum_{\alpha\in I} (e_\alpha \otimes 1) \otimes (1 \otimes e^\alpha).
\end{equation} 
\end{defn}

\begin{prop}
The universal $R$-matrix satisfies the {\normalfont{Yang-Baxter equation}}
\begin{equation}\label{yang_baxter_equation} 
 R_{12} R_{13} R_{23} = R_{23} R_{13} R_{12} .
\end{equation}
where we use a notation for which $R_{12} = R \otimes (1\otimes 1)$, $R_{23} = (1\otimes 1)\otimes R$ and $R_{13} = \sum_{\alpha\in I} e_\alpha \otimes (1\otimes 1) \otimes e^\alpha$. 
\end{prop}

The Drinfeld double $D(\cA)$ can be related to the Heisenberg algebras in terms of an algebra morphism. Lets consider again the Heisenberg double $H(\cA)$, as described in section \ref{chapter2}. Let us recall that one can regard it as an algebra of elements $\{e^\beta \otimes e_\alpha\}_{\alpha,\beta\in I}$ subjected to the set of relations
\begin{align}
 (e^\alpha \otimes e_\beta)(e^\gamma \otimes e_\delta) &= \sum_{\epsilon,\pi,\rho,\sigma,\tau\in I} (-1)^{|\beta||\gamma| + |\pi||\epsilon| + |\pi||\alpha| + |\epsilon|} m^\gamma_{\pi\epsilon} m^\tau_{\rho\delta} \mu^{\epsilon\rho}_\beta \mu^{\alpha\pi}_\sigma e^\sigma \otimes e_\tau.
\end{align}
Moreover, we can also construct an additional Heisenberg double $H(\cA^*) = \{{\tilde e}_\alpha \otimes {\tilde e}^\beta \}_{\alpha,\beta\in I}$ starting from the dual algebra $\cA^*$. It has the following relations
\begin{equation}\begin{split}
({\tilde e}_\alpha \otimes \tilde e^\beta)(\tilde e_\gamma \otimes {\tilde e}^\delta) &= \sum_{\epsilon,\pi,\rho,\sigma,\tau\in I} (-1)^{|\rho||\pi|+|\rho||\epsilon|+|\pi||\delta|} \mu^{\rho\epsilon}_\gamma \mu^{\pi\delta}_\tau m^\beta_{\epsilon\pi} m^\sigma_{\alpha\rho} \tilde e_\sigma \otimes \tilde e^\tau .
\end{split}\end{equation}
We will denote the flipped (i.e. the one with the tensor factors reversed) canonical element of this Heisenberg double as ${\tilde S} = {\tilde e}_\alpha \otimes {\tilde e}^\alpha$. It satisfies a ``reversed'' pentagon equation of the form
\begin{align}
 {\tilde S}_{12} {\tilde S}_{23} = {\tilde S}_{23} {\tilde S}_{13} {\tilde S}_{12} .
\end{align}
We can relate those 2 algebras by the means of the following proposition:
\begin{prop}
 There exists an algebra anti-isomorphism $\xi: H(\cA^*) \to H(\cA)$ given by
 \begin{align}
  \xi({\tilde e}_\alpha) = (-1)^{c|\alpha|} \gamma_\alpha^\beta e_\beta, && \xi({\tilde e}^\alpha) = (-1)^{(c+1)|\alpha|} (\gamma^{-1})^\alpha_\beta e^\beta,
 \end{align}
 where $c=0,1$.
\end{prop}
The anti-isomorphism can be implemented on representation spaces in terms of super-transposition (i.e. the graded analogue of ordinary transposition). The super-transposition for square even $(n|m)$-matrices, i.e. linear transformations belonging to the space of $Hom(C^{(n|m)},C^{(n|m)})$, is given by 
\begin{equation}
 \left( \begin{array}{cc} A & B\\C& D \end{array} \right)^\text{st} = \left( \begin{array}{cc} A^\text{t} & C^\text{t}\\-B^\text{t}& D^\text{t} \end{array} \right) ,
\end{equation}
where $A\in Hom(\mathbb{C}^n,\mathbb{C}^n), D\in Hom(\mathbb{C}^m,\mathbb{C}^m)$ are even, and $B\in Hom(\mathbb{C}^m,\mathbb{C}^n), C\in Hom(\mathbb{C}^n,\mathbb{C}^m)$ are odd, and $^\text{t}$ denotes ordinary, not graded matrix transposition. Moreover, the transposition on $L^2(\mathbb{R})$ is implemented by the following action on the momentum and position operators: $p^\text{t} = -p, q^\text{t} = q$. \\

\noindent Then, we claim that the following proposition is true: 
\begin{prop} A map $\eta: D(\cA) \to H(\cA)\otimes H(\cA^*)$ defined as follows \label{drinfeld_heisenberg_homomorphism}
\begin{align}
 &\eta^{(a,b)}(e_\alpha\otimes1) = (-1)^{a |\beta| + b|\gamma|} \mu_\alpha^{\beta\gamma} e_\beta\otimes\tilde e_\gamma, && \eta^{(a,b)}(1\otimes e^\alpha) = (-1)^{a' |\beta| + b'|\gamma|} m^\alpha_{\gamma\beta} e^\beta\otimes\tilde e^\gamma, 
\end{align}
for the choice of the parameters $a,a',b,b'\in\mathbb{Z}_{\geq0}$ such that
\begin{align}
 (-1)^{a+a'}=-1, && (-1)^{b+b'} = 1 ,
\end{align}
is an algebra homomorphism.
\end{prop}
This fact can be checked by a direct calculation. \\

Using the morphism $\eta$ from the proposition \ref{drinfeld_heisenberg_homomorphism} in addition to the canonical element $S$ for the Heisenberg double $H(\cA)$ and the (flipped) canonical element ${\tilde S}$ for $H(\cA^*)$, one can define the following two elements $S' = (-1)^{(a+b'+1)|\alpha|} {\tilde e}_\alpha \otimes e^\alpha$ and $S'' = (-1)^{(a+b') |\alpha|} { e}_\alpha \otimes {\tilde e}^\alpha$.
It can be shown that they satisfy a set of 6 pentagon-like equations
\begin{align}
 S'_{12} S'_{13} S_{23} = S_{23} S'_{12}, && {\tilde S}_{12} S'_{23} = S'_{23} S'_{13} {\tilde S}_{12}, \nonumber \\
 S_{12} S''_{13} S''_{23} = S''_{23} S_{12}, && S''_{12} {\tilde S}_{23} = {\tilde S}_{23} S''_{13} S''_{12}, \\
 S'_{12} {\tilde S}_{13} S''_{23} = S''_{23} S'_{12}, && S''_{12} S'_{23} = S'_{23} S_{13} S''_{12}. \nonumber
\end{align}
Then, we claim that one can construct the R-matrix of the Drinfeld double $D(\cA)$ as follows
\begin{prop}
 Under an algebra map $\eta$ one has the following relation
 \begin{equation}
  (\eta^{(a,b)}\otimes \eta^{(a,b)}) R = S''_{14} S_{13} {\tilde S}_{24} S'_{23} .
 \end{equation}
\end{prop}

\begin{re}
 To keep the notation compact, from now on we will denote the elements $1\otimes e^\alpha $ and $e_\alpha\otimes1 $ of the Drinfeld double $D(\cA)$ simply as $e^\alpha$ and $e_\alpha$ respectively.
\end{re}

\begin{exmp} Drinfeld double of the Borel half of $U_q(sl(2))$ \end{exmp} \vspace{-0.25cm}

Using the definitions above, one can obtain the Drinfeld double $D(\cA)$ commutation relations for $\cA=\mathcal{B}(U_q(sl(2)))$
\begin{align} 
\begin{aligned}
 &[H,E] = -ibE, &&\qquad [H,F] = +ibF,\\
 &[\hat H, E] = -ibE, &&\qquad [\hat H,F] = +ibF,\\
 &[H,\hat H] = 0, &&\qquad [E,F] = (q-q^{-1})(e^{2\pi b H}-e^{-2\pi b {\hat H}}) ,
 \end{aligned}
\end{align}
with the coproduct
\begin{equation}
 \begin{aligned}
 &\Delta(H) = H\otimes1 + 1\otimes H, && \Delta(E) = E\otimes e^{2\pi b H} + 1\otimes E, \\
 &{ \Delta}(\hat H) = \hat H\otimes1 + 1\otimes \hat H,
 &&{ \Delta}(F) = F\otimes e^{-2\pi b {\hat H}} + 1\otimes F,
\end{aligned}
\end{equation}
and the antipode
\begin{align}
\begin{aligned}
 \gamma(H) &= -H, && \gamma(E) = q E e^{-2\pi b H} ,\\
 \gamma(\hat H) &= - \hat H, && \gamma(F) = q F e^{2\pi b \hat H}. 
 \end{aligned}
\end{align}
Using the map $\eta$ as well as the algebra anti-homomorphism $\xi$, one can obtain a representation $\pi_{D(\cA)}: D(\cA) \to L^2(\mathbb{R})^{\otimes 2}$ from the representation \eqref{representation_sl2}
\begin{align} 
\begin{aligned}
 &\pi_{D(\cA)}(H) = \mathsf{p}_1+\mathsf{p}_2, && \pi_{D(\cA)}(E) = e^{2\pi b (\mathsf{p}_2+\mathsf{x}_2)} + e^{2\pi b (\mathsf{x}_1+\mathsf{p}_2)},\\
 &\pi_{D(\cA)}({\hat H}) = \mathsf{x}_1-\mathsf{x}_2, && \pi_{D(\cA)}(F) = e^{-2\pi b (\mathsf{x}_1+\mathsf{p}_2)} + e^{2\pi b (\mathsf{p}_1-\mathsf{x}_1)} ,
\end{aligned} 
\end{align}
where $\mathsf{p},\mathsf{x}$ are the momentum and position operators satisfying $[\mathsf{p},\mathsf{x}]=\frac{1}{2\pi i}$. Moreover, the universal $R$-matrix realised using this representation is given by
\begin{equation}\begin{aligned}
 (\pi_{D(\cA)} \otimes \pi_{D(\cA)}) R &= e^{2\pi i(\mathsf{p}_1+\mathsf{p}_2)(\mathsf{x}_3-\mathsf{x}_4)} g_\ub^{-1}(e^{2\pi b(\mathsf{x}_1+\mathsf{p}_2-\mathsf{x}_3-\mathsf{p}_4)}) g_\ub^{-1}(e^{2\pi b(\mathsf{x}_1+\mathsf{p}_2+\mathsf{p}_3-\mathsf{x}_3)}) \times \\
 & \qquad \times g_\ub^{-1}(e^{2\pi b(\mathsf{p}_2+\mathsf{x}_2-\mathsf{x}_3-\mathsf{p}_4)}) g_\ub^{-1}(e^{2\pi b(\mathsf{p}_2+\mathsf{x}_2+\mathsf{p}_3-\mathsf{x}_3)}).
\end{aligned}\end{equation}

\begin{exmp} Drinfeld double of the Borel half of $U_q(osp(1|2))$ \end{exmp} \vspace{-0.25cm}

As another instructive example, one can consider the Drinfeld double $D(\cA)$ for $\cA=\mathcal{B}(U_q(osp(1|2)))$. Its commutation relations are as follows
\begin{align} 
\begin{aligned}
 &[H,\vp] = -ib\vp, &&\qquad [H,\vm] = +ib\vm,\\
 &[\hat H, \vp] = -ib\vp, &&\qquad [\hat H,\vm] = +ib\vm,\\
 &[H,\hat H] = 0, &&\qquad \{\vp,\vm\} = (q^\frac{1}{2}+q^{-\frac{1}{2}})(e^{\pi b H}-e^{-\pi b {\hat H}}) .
 \end{aligned}
\end{align}
with the coproduct
\begin{equation}
 \begin{aligned}
 &\Delta(H) = H\otimes1 + 1\otimes H, && \Delta(\vp) = \vp\otimes e^{\pi b H} + 1\otimes \vp, \\
 &{ \Delta}(\hat H) = \hat H\otimes1 + 1\otimes \hat H,
 &&{ \Delta}(\vm) = \vm\otimes e^{-\pi b {\hat H}} + 1\otimes \vm,
\end{aligned}
\end{equation}
and the antipode
\begin{align}
\begin{aligned}
 \gamma(H) &= -H, && \gamma(\vp) = q^\frac{1}{2} \vp e^{-\pi b H} ,\\
 \gamma(\hat H) &= - \hat H, && \gamma(\vm) = q^{-\frac{1}{2}} \vm e^{\pi b \hat H}. 
 \end{aligned}
\end{align}
Using the map $\eta$ as well as the algebra anti-homomorphism $\xi$, one can obtain a representation $\pi_{D(\cA)}: D(\cA) \to L^2(\mathbb{R})^{\otimes 2} \otimes (\mathbb{C}^{1|1})^{\otimes 2}$ from the representation \eqref{representation_osp12}
\begin{align} 
\begin{aligned}
 &\pi_{D(\cA)}^{(u,v)}(H) = (\mathsf{p}_1+\mathsf{p}_2) \mathbb{I}_2, \\ 
 & \pi_{D(\cA)}^{(u,v)}(\vp) = (-1)^u i \left[ e^{\pi b(\mathsf{x}_1+\mathsf{p}_2)} \left( \begin{array}{cc} 0 & 1\\1& 0 \end{array} \right) \otimes \mathbb{I}_2 + (-1)^{v} e^{\pi b(\mathsf{p}_2+\mathsf{x}_2)} \mathbb{I}_2 \otimes \left( \begin{array}{cc} 0 & 1\\-1& 0 \end{array} \right) \right],\\
 &\pi_{D(\cA)}^{(u,v)}({\hat H}) = (\mathsf{x}_1-\mathsf{x}_2) \mathbb{I}_2, \\
 & \pi_{D(\cA)}^{(u,v)}(\vm) = (-1)^{u+1} i \left[ e^{\pi b (\mathsf{p}_1-\mathsf{x}_1)} \left( \begin{array}{cc} 0 & 1\\1& 0 \end{array} \right) \otimes \mathbb{I}_2 + (-1)^{v} e^{-\pi b (\mathsf{x}_1+\mathsf{p}_2)} \mathbb{I}_2 \otimes \left( \begin{array}{cc} 0 & 1\\-1& 0 \end{array} \right) \right],
\end{aligned} 
\end{align}
where $u,v=0,1$ are free parameters, $\mathsf{p},\mathsf{x}$ are the momentum and position operators satisfying $[\mathsf{p},\mathsf{x}]=\frac{1}{\pi i}$ and the $(1|1)$-dimensional identity matrix is denoted as $\mathbb{I}_2 = \left( \begin{array}{cc} 1 & 0\\0& 1 \end{array} \right)$. Moreover, the previously not obtained universal $R$-matrix realised using this representation is given by
\begin{small}\begin{equation}\begin{aligned}
 & \left(\pi_{D(\cA)}^{(u,v)} \otimes \pi_{D(\cA)}^{(u,v)} \right) R = \frac{1}{16} e^{\pi i(\mathsf{p}_1+\mathsf{p}_2)(\mathsf{x}_3-\mathsf{x}_4)} \times \\
 &\times \left[ h_+(\mathsf{x}_1+\mathsf{p}_2-\mathsf{x}_3-\mathsf{p}_4)) \mathbb{I}_2^{\otimes4} + (-1)^{v} i h_-(\mathsf{x}_1+\mathsf{p}_2-\mathsf{x}_3-\mathsf{p}_4) \left( \begin{array}{cc} 0 & 1\\1& 0 \end{array} \right) \otimes \mathbb{I}_2 \otimes \mathbb{I}_2 \otimes \left( \begin{array}{cc} 0 & 1\\-1& 0 \end{array} \right) \right] \times \\
 &\times \left[ h_+(\mathsf{x}_1+\mathsf{p}_2+\mathsf{p}_3-\mathsf{x}_3)) \mathbb{I}_2^{\otimes4} - i h_-(\mathsf{x}_1+\mathsf{p}_2+\mathsf{p}_3-\mathsf{x}_3) \left( \begin{array}{cc} 0 & 1\\1& 0 \end{array} \right) \otimes \mathbb{I}_2 \otimes \left( \begin{array}{cc} 0 & 1\\1& 0 \end{array} \right) \otimes \mathbb{I}_2 \right] \times \\
 &\times \left[ h_+(\mathsf{p}_2+\mathsf{x}_2-\mathsf{x}_3-\mathsf{p}_4)) \mathbb{I}_2^{\otimes4} + i h_-(\mathsf{p}_2+\mathsf{x}_2-\mathsf{x}_3-\mathsf{p}_4) \, \mathbb{I}_2 \otimes \left( \begin{array}{cc} 0 & 1\\-1& 0 \end{array} \right) \otimes \mathbb{I}_2 \otimes \left( \begin{array}{cc} 0 & 1\\-1& 0 \end{array} \right) \right] \times \\
 &\times \left[ h_+(\mathsf{p}_2+\mathsf{x}_2+\mathsf{p}_3-\mathsf{x}_3)) \mathbb{I}_2^{\otimes4} + (-1)^{v} i h_-(\mathsf{p}_2+\mathsf{x}_2+\mathsf{p}_3-\mathsf{x}_3) \, \mathbb{I}_2 \otimes \left( \begin{array}{cc} 0 & 1\\-1& 0 \end{array} \right) \otimes \left( \begin{array}{cc} 0 & 1\\1& 0 \end{array} \right) \otimes \mathbb{I}_2  \right] ,
\end{aligned}\end{equation}
\end{small}
where $h_\pm(x) = e_\R^{-1}(x) \pm e_\NS^{-1}(x)$. 

\section{Outlook} \label{chapter6}

Let us conclude our paper by mentioning some interesting directions for future work suggested by our results. 

As we mentioned in the introduction, a 1-parameter family of infinite dimensional representations $\mathcal{P}_\alpha, \, \alpha \in \frac{1}{2}(b+b^{-1})+i\mathbb{R}$ of $U_q(sl(2)), \, q=e^{i\pi b^2}$ has been studied by Bytsko, Ponsot, and Teschner \cite{PT1,PT2, Bytsko:2002br} in connection to the Liouville theory. $\mathcal{P}_\alpha$ represents $U_q(sl(2))$ on a space of analytic function defined on a strip around $\{x \in \mathbb{C}: |\text{Im}(x)| < \frac{b}{2}\}$ which possess a Fourier transform that is meromorphic on $\mathbb{C}$ with a specified set of allowed poles. The generators $K,E,F$ of the quantum group are realised as positive, self-adjoint operators. The family $\mathcal{P}_\alpha$ is closed under the the tensor product and the calculation of the 3j- and 6j-symbols allowed to make a connection with the fusion matrices of the Liouville theory. Moreover, on the same representation space acts a representation of $U_{\tilde q}(sl(2))$ for ${\tilde q}=e^{i\pi b^2}$, and therefore $\mathcal{P}_\alpha$ constitute representations of the modular double of the quantum group introduced in a sense of Faddeev \cite{Faddeev:1999fe}. This modularity property, which ensures a self-duality of the exchange $b \to \frac{1}{b}$, is crucial for the interpretation in terms of Liouville theory, which exhibits the same symmetry. The representations $\mathcal{P}_\alpha$ were also found in during the study of the spectral problem of Dehn twists in quantum Teichm\"uller theory utilising an algebra map similar to the one in proposition \ref{drinfeld_heisenberg_homomorphism} --- it differed form it however by a twisting of the co-product. 

This results for the class of infinite-dimensional representations $\mathcal{P}_\alpha$ of the quantum group $U_q(sl_2)$ has been generalised in \cite{gus, gus-cluster} for the case of the higher rank $U_q(sl_{n+1})$ quantum groups using the cluster algebras methods. They describe the algebraic ingredients of a proof of the conjecture of Frenkel and Ip \cite{Ip2} that the category of the representations $\mathcal{P}_\lambda$ of the quantum group $U_q(sl_{n+1})$ is closed under tensor products.

In the context of $N=1$ supersymmetric Liouville theory, the attempt to find an $U_q(osp(1|2))$ analogue of the representations $\mathcal{P}_\alpha$ resulted obtaining the 6j-symbols which reproduce only the Neveu-Schwarz sector of the theory \cite{Hadasz:2013bwa,Pawelkiewicz:2013wga}. Given our results regarding the Heisenberg double of the Borel half of $U_q(osp(1|2))$, we conjecture that the proper analogue of the representations $\mathcal{P}_\alpha$, i.e. the one which will encode the entirety of the structure of $N=1$ supersymmetric Liouville theory, is a 1-parameter family of infinite-dimensional representations of $U_q(osp(1|2))$ modeled on a dense subspace of the Hilbert space $L^2(\mathbb{R}) \otimes (\mathbb{C}^{1|1})^{\otimes2}$ and given by
\begin{equation}\begin{aligned}
& \pi_\alpha^{(u,v)}(K) = e^{\pi b \mathsf{p}} \mathbb{I}_2 \otimes \mathbb{I}_2, \\
& \pi_\alpha^{(u,v)}(v ^{(\pm)}) =  (-1)^u i q^{\pm\frac{1}{4}} \left[ e^{\frac{\pi b}{2}( \mp 2 \mathsf{x} + \mathsf{p} \pm \alpha)} \left( \begin{array}{cc} 0 & 1\\1& 0 \end{array} \right) \otimes \mathbb{I}_2 \pm (-1)^v e^{\frac{\pi b}{2}( \mp 2 \mathsf{x} - \mathsf{p} \mp \alpha)} \mathbb{I}_2 \otimes \left( \begin{array}{cc} 0 & 1\\1& 0 \end{array} \right) \right],
\end{aligned}\end{equation}
where $u,v=0,1$ are free parameters, $q=\exp[i\pi b^2]$ and where the momentum and position operators $\mathsf{p},\mathsf{x}$ acting on $L^2(\mathbb{R})$ satisfy $[\mathsf{p},\mathsf{x}]=\frac{1}{\pi i}$. $\pi_\alpha$ is a representation of an $U_q(osp(1|2))$ algebra defined by the (anti-)commutation relations
\begin{align} 
\begin{aligned}
 &K v^{(\pm)} = q^{\pm1} v^{(\pm)} K, \qquad \{\vp,\vm\} = (q^\frac{1}{2}+q^{-\frac{1}{2}})(K-K^{-1}) ,
 \end{aligned}
\end{align}
and the co-product
\begin{equation}
 \begin{aligned}
 &\Delta(K) = K\otimes K, && \Delta(v^{(\pm)}) = v^{(\pm)} \otimes K + K^{-1}\otimes v^{(\pm)} .
\end{aligned}
\end{equation}
The novel $R$-matrix for this family of representations is given by
\begin{equation}
 (\pi_\alpha^{(u,v)}\otimes\pi_\beta^{(u,v)}) R = \frac{1}{16} F^{-1} e^{i\pi (\mathsf{p}_1+\alpha)(\mathsf{p}_2-\beta)} S_1 S_2 S_3 S_4 F^{-1} ,
\end{equation}
where
\begin{align*}
 F &= \exp\left\{ \frac{i\pi}{4} [(\mathsf{p}_1+\alpha)(\mathsf{p}_2-\beta) - (\mathsf{p}_1-\alpha)(\mathsf{p}_2+\beta)] \right\} , \\
 S_1 &= h_+(\mathsf{p}_1+\mathsf{x}_1+\mathsf{p}_2+\mathsf{x}_2)) \mathbb{I}_2^{\otimes4} - (-1)^{v} i h_-(\mathsf{p}_1+\mathsf{x}_1+\mathsf{p}_2+\mathsf{x}_2) \left( \begin{array}{cc} 0 & 1\\1& 0 \end{array} \right) \otimes \mathbb{I}_2 \otimes \mathbb{I}_2 \otimes \left( \begin{array}{cc} 0 & 1\\1& 0 \end{array} \right), \\
 S_2 &= h_+(\mathsf{p}_1+\mathsf{x}_1+\mathsf{p}_2+\beta)) \mathbb{I}_2^{\otimes4} - i h_-(\mathsf{p}_1+\mathsf{x}_1+\mathsf{p}_2+\beta) \left( \begin{array}{cc} 0 & 1\\1& 0 \end{array} \right) \otimes \mathbb{I}_2 \otimes \left( \begin{array}{cc} 0 & 1\\1& 0 \end{array} \right) \otimes \mathbb{I}_2 ,\\
 S_3 &= h_+(\mathsf{x}_1+\alpha+\mathsf{p}_2+\mathsf{x}_2)) \mathbb{I}_2^{\otimes4} + i h_-(\mathsf{x}_1+\alpha+\mathsf{p}_2+\mathsf{x}_2) \mathbb{I}_2 \otimes \left( \begin{array}{cc} 0 & 1\\1& 0 \end{array} \right) \otimes \mathbb{I}_2 \otimes \left( \begin{array}{cc} 0 & 1\\1& 0 \end{array} \right), \\
 S_4 &= h_+(\mathsf{x}_1+\alpha+\mathsf{p}_2+\beta)) \mathbb{I}_2^{\otimes4} + (-1)^{v} i h_-(\mathsf{x}_1+\alpha+\mathsf{p}_2+\beta) \mathbb{I}_2 \otimes \left( \begin{array}{cc} 0 & 1\\1& 0 \end{array} \right) \otimes \left( \begin{array}{cc} 0 & 1\\1& 0 \end{array} \right)  \otimes \mathbb{I}_2.
\end{align*}
It would be interesting to perform a harmonic analysis for this family of representations and calculate the 3j- and 6j-symbols. In particular, of interest is whether 6j-symbols could reproduce the entire fusion matrix of the $N=1$ supersymmetric Liouville theory.

\newpage


\acknowledgments

We are very grateful to J\"org Teschner for explanations, suggestions and many helpful discussions and comments. We also thank Rinat Kashaev for stimulating discussions.

This project was initiated when NA and MP was supported by the German Science Foundation (DFG) within the Research Training Group 1670 ''Mathematics Inspired by String Theory and QFT''. The work of N.A. was supported by Max Planck Institute of Mathematics (MPIM)in Bonn. The work of M.P. was supported by the European Research Council (advanced grant NuQFT). We also thank and acknowledges the AEC centre at University of Bern and IHES for their hospitality during this project. 
\\


\appendix

\section{Special functions}\label{appendixA}

Quantum dilogarithm plays a key role in the constructions described in this paper. In this appendix we review the Faddeev's quantum dilogarithm and its most important properties. We collected the different definitions of related special functions which one may face in the references.
\subsection{Faddeev's quantum dilogarithm}\label{Appendix.A.1}

The basic special function that appears in the context of the infinite dimensional representations of the Heisenberg double of the quantum plane is Barnes' double gamma function \cite{brane}. The double gamma function is defined as
\begin{align*}
\log\Gamma_2(z|\omega):= \frac
{\partial}{\partial s} \left( \sum_{m_1,m_2\in \bZ_{\geq0}}{(z+m_1\omega_1+m_2\omega_2)}^{-s} \right)_{s=0},
\end{align*}
and using which one can define
\begin{align*}
\Gamma_\ub(x):=\Gamma_2(x|\ub,\ub^{-1}).
\end{align*}
For $\mathfrak{Re}x>0$ it admits an integral representation,
\begin{eqnarray}\label{def_specfunc}
\log \Gamma_\ub(x) = \int_0^\infty \frac{\ud t}{t}
\left[ \frac{e^{-x t} - e^{-\frac{Q}{2} t}}{(1 - e^{-t \ub})
	(1 - e^{-\frac{t}{\ub}})} - \frac{\left( \frac{Q}{2}-x\right)^2}{2 e^t}
- \frac{\frac{Q}{2} - x}{t} \right] ,
\end{eqnarray}
where $Q = \ub + \frac{1}{\ub}$. One can analytically continue $\Gamma_\ub$
to a meromorphic function defined on the entire complex plane
$\mathbb{C}$. The most important property of $\Gamma_\ub$ is its
behavior with respect to shifts by $\ub^\pm$,
\begin{equation}
\Gamma_\ub(x+\ub) = \frac{\sqrt{2\pi} \ub^{\ub x-\frac{1}{2}}}{\Gamma_\ub(bx)}
\Gamma_\ub(x)\quad , \quad
\Gamma_\ub(x+\ub^{-1}) = \frac{\sqrt{2\pi} \ub^{-\frac{\ub}{x}+\frac{1}{2}}}
{\Gamma_\ub(\frac{x}{\ub})} \Gamma_\ub(x)\ .
\end{equation}
These shift equation allows us to calculate residues of the poles of
$\Gamma_\ub$. When $x\to0$, for instance, one finds
\begin{equation}\label{residue}
\Gamma_\ub(x)= \frac{\Gamma_\ub(Q)}{2\pi x} + O(1).
\end{equation}
From Barnes' Double Gamma function we can build  other important
special functions,
\begin{align}
\Upsilon_\ub(x)&:=\frac{1}{\Gamma_\ub(x)\Gamma_\ub(Q-x)},\\[2mm]
S_\ub(x) &:= \frac{\Gamma_\ub(x)}{\Gamma_\ub(Q-x)},\\[2mm]
G_\ub(x) &:= e^{-\frac{i\pi}{2} x(Q-x)} S_\ub(x), \label{GG}\\[2mm]
w_\ub(x) &:=e^{\frac{\pi i}{2}(\frac{Q^2}{4}+x^2)}G_\ub(\frac{Q}{2}-ix),\\[2mm]
g_{\ub}(x)&:=\frac{\zeta_\ub}{G_\ub(\frac{Q}{2}+\frac{1}{2\pi i \ub}\log x)},\\[2mm]
e_{\ub}(x)&:=\frac{\zeta_\ub}{G_\ub(\frac{Q}{2}-i x)},
\end{align}
where $\zeta_\ub = \exp[-\frac{i\pi}{4}-\frac{i\pi}{12}(\ub^2+\ub^{-2})]$
One refers to the function $S_\ub$ as double sine function. It is defined by the following integral representation
\begin{equation}
\log S_\ub(x)=\int_0^{\infty} \frac{\ud t}{it} \left(\frac {\sin{2xt}}{2\sinh{\ub t}\sinh{\ub^{-1}t}-\frac{x}{t}} \right).
\end{equation}
The $S_\ub$ function is meromorphic with poles and zeros in
\begin{align*}
S_\ub(x) = 0 &\Leftrightarrow x = Q + n \ub + m \ub^{-1}, &&\quad n,m \in \mathbb{Z}_{\geq0}\ , \\[2mm]
S_\ub(x)^{-1} = 0 &\Leftrightarrow x = -n \ub -m \ub^{-1}, &&\quad n,m \in \mathbb{Z}_{\geq0}\ .
\end{align*}
Other most important properties are as follows:
\begin{align}
&\text{Functional equation(Shift):} && S_\ub(x-i\ub/2)=2\cosh {(\pi \ub x)}S_\ub(x+i\ub/2)\\
&\text{Self-duality:} &&S_\ub(x)=S_{1/\ub}(x)\\
&\text{Inversion relation(Reflection):}&& S_\ub(x)S_\ub(-x)=1\\
&\text{Unitarity:} && \overline{S_\ub(x)}=1/S_\ub(\overline{x})\\
&\text{Residue:} && \text{res}_{x=iQ/2}S_\ub(x)= e^{-\frac{i\pi}{12}(1+{Q}^2)}{(2\pi i)}^{-1}.
\end{align}
From the relation between the special functions and the shift property of Barnes' double Gamma
function it is easy to derive the following shift and reflection
properties of $G_\ub$,
\begin{align}
&G_\ub(x+\ub) = (1-e^{2 \pi i \ub x}) G_\ub(x)\ ,\label{shift} \\[2mm]
&G_\ub(x) G_\ub(Q-x) = e^{\pi i x(x-Q)}\ .\label{reflection}
\end{align}
We also need to the asymptotic behavior of the function $G_\ub$ along
the imaginary axis,
\begin{equation}\begin{split}\label{asymptotic}
G_\ub(x) \ \sim &\  \zeta_\ub \ , \qquad \qquad \ \ \ \  \mathfrak{Im} x \to +\infty,\\[2mm]
G_\ub(x) \sim & \ \zeta_\ub^{-1} e^{i\pi x(x-Q)}\, , \quad \mathfrak{Im}x \to -\infty.
\end{split}
\end{equation}

The Fadeev's quantum dilogarithm function in addition to the relation with the $G_\ub$ function has the following integral representation
\begin{equation}
 g_\ub \left(\frac{1}{2\pi \ub} \log x \right) = \exp \left[\int_{\mathbb{R}+i0} \ud w \frac{e^{-2ixw}}{4 \sinh(w\ub)\sinh(w\slash \ub)}\right],
\end{equation}
and
\begin{align*}
& g_\ub(e^{2\pi \ub r}) = \int\ud t\, e^{2\pi i t r}\frac{e^{-i\pi t^2}}{G_\ub(Q+it)}, \\
& g_\ub^{-1}(e^{2\pi \ub r}) = \int\ud t\, e^{2\pi i t r}\frac{e^{-\pi tQ}}{G_\ub(Q+it)},
\end{align*}
The shift and reflection relations that it satisfies are as follows
\begin{align*}
& g_\ub( e^{-i\pi \ub^2} x) = (1+x) g_\ub( e^{+i\pi \ub^2} x) ,\\
& g_\ub(e^{2\pi i\ub x}) g_\ub(e^{-2\pi i\ub x}) = e^{\frac{i\pi Q^2}{4}} \zeta_\ub^2 e^{i\pi x^2}.
\end{align*}
Also, for non-commutative variables $U,V$ such that $UV = q^2 VU$ where $q=e^{i\pi \ub^2}$ it satisfies the pentagon relation
\begin{equation}\label{pentagonfaddeev_gb_function}
g_\ub(U) g_\ub(V) = g_\ub(V) g_\ub(q^{-1} UV) g_\ub(U).
\end{equation}
The pentagon equation can be equivalently expressed as the Ramanujan summation formula \cite{PT2,FKV,Vol}
\begin{equation} \label{Ramanujan}
\int_{-i\infty}^{i\infty} \frac{\ud\tau}{i} e^{2 \pi i \tau \beta} \frac{G_\ub(\tau + \alpha)}{G_\ub(\tau + Q)} = \frac{G_\ub(\alpha) G_\ub(\beta)}{G_\ub(\alpha+\beta)}.
\end{equation}
Moreover, for the function $e_\ub$ we have the following shift and reflection relations
\begin{align}
e_\ub\left(x-\frac{i \ub^{\pm1}}{2}\right) &= (1+e^{2\pi \ub^{\pm1}x})e_\ub\left(x+\frac{i\ub^{\pm1}}{2}\right), \\
e_\ub(x)e_\ub(-x) &= e^{-i\pi(1-Q^2/2)\slash 6} e^{i\pi x^2}.
\end{align}
The asymptotic behaviour of the function $e_\ub$ along the real axis 
\begin{equation} \label{asymptotic}
e_\ub(x) = \left\{ \begin{array}{ll}
1 &, x\to -\infty\\
e^{-i\pi(1-Q^2/2)\slash 6} e^{i\pi x^2} & ,x\to +\infty\\
\end{array} \right.
\end{equation}
Also, we know that for self-adjoint operators $\mathsf{p},\mathsf{x}$ such that $[\mathsf{p},\mathsf{x}]=\frac{1}{2\pi i}$  we have the following variant of the pentagon relation 
\begin{equation}\label{pentagonfaddeev}
e_\ub(\mathsf{p}) e_\ub(\mathsf{x}) = e_\ub(\mathsf{x}) e_\ub(\mathsf{x}+\mathsf{p}) e_\ub(\mathsf{p}) .
\end{equation}

\subsection{Supersymmetric non-compact quantum dilogarithm}\label{Appendix.A.2}

Now, we will consider special functions related to the supersymmetric analogue of the Faddeev's quantum dilogarithm. We can define the supersymmetric analogues of double gamma functions
\begin{align*}
\Gamma_\NS(x) &= \Gamma_\ub\left(\frac{x}{2}\right) \Gamma_\ub
\left(\frac{x+Q}{2}\right),\\[2mm]
\Gamma_\R(x) &= \Gamma_\ub\left(\frac{x+\ub}{2}\right) \Gamma_\ub
\left(\frac{x+\ub^{-1}}{2}\right).
\end{align*}
Furthermore, let us define
\begin{equation}
\begin{array}{rlrl}
S_\NS(x) = \frac{\Gamma_\NS(x)}{\Gamma_\NS(Q-x)} , \quad &
G_\NS(x) = \zeta_0 e^{-\frac{i\pi}{4} x(Q-x)} S_\NS(x), \\[2mm]
S_\R(x) = \frac{\Gamma_\R(x)}{\Gamma_\R(Q-x)}, \quad &
G_\R(x) = e^{-\frac{i\pi}{4}} \zeta_0 e^{-\frac{i\pi}{4}
	x(Q-x)} S_\R(x),
\end{array}
\end{equation}
where $\zeta_0 = \exp(-i\pi Q^2/8)$. As for $S_\ub$, the functions $S_\R(x)$ and $S_\NS(x)$ are meromorphic with poles and zeros in
\begin{align*}\
S_\R(x) = 0 &\Leftrightarrow x = Q + n \ub + m \ub^{-1}, \quad n,m \in \mathbb{Z}_{\geq0}, m+n\in2\mathbb{Z}+1, \\
S_\NS(x) = 0 &\Leftrightarrow x = Q + n \ub + m \ub^{-1}, \quad n,m \in \mathbb{Z}_{\geq0}, m+n\in2\mathbb{Z},\\
S_\R(x)^{-1} = 0 &\Leftrightarrow x = -n \ub -m \ub^{-1}, \quad n,m \in \mathbb{Z}_{\geq0}, m+n\in2\mathbb{Z}+1, \\
S_\NS(x)^{-1} = 0 &\Leftrightarrow x = -n \ub -m \ub^{-1}, \quad n,m \in \mathbb{Z}_{\geq0}, m+n\in2\mathbb{Z}.
\end{align*}
We state the shift and reflection properties of the functions $G_\NS$ and $G_\R$
\begin{align}
G_\R(x+\ub^{\pm1}) &= (1- e^{ \pi i \ub^{\pm1} x}) G_{\NS}(x),
\label{shift_super1}\\[2mm]
G_\NS(x+\ub^{\pm1}) &= (1+ e^{ \pi i \ub^{\pm1} x}) G_{\R}(x),
\label{shift_super1}\\[2mm]
\label{reflection_super2}
G_\R(x) G_\R(Q-x) &= e^{-\frac{i\pi}{2}} \zeta_0^2
e^{\frac{\pi i}{2} x(x-Q)}\ , \\
\label{reflection_super2}
G_\NS(x) G_\NS(Q-x) &= \zeta_0^2
e^{\frac{\pi i}{2} x(x-Q)}\ .
\end{align}
We define the supersymmetric analogues of Faddeev's quantum dilogarithm function as
\begin{align}\begin{aligned}
e_\R(x) &= e_\ub\left(\frac{x}{2}+\frac{i }{4}(\ub-\ub^{-1})\right)e_\ub\left(\frac{x}{2}-\frac{i }{4}(\ub-\ub^{-1})\right),\\
e_\NS(x) &= e_\ub\left(\frac{x}{2}+\frac{i }{4}(\ub+\ub^{-1})\right)e_\ub\left(\frac{x}{2}-\frac{i }{4}(\ub+\ub^{-1})\right) .
\end{aligned}\end{align}
and relate them to the double sine function in a way as follows
\begin{equation}\begin{aligned}
e_\R(x) &=\ \frac{\zeta_\ub^2}{G_\R(-ix + \frac{Q}{2})}, \\
e_\NS(x) &=\ \frac{\zeta_\ub^2}{ G_\NS(-ix + \frac{Q}{2})}.
\end{aligned}\end{equation}
In addition, the functions $e_\R$ and $e_\NS$ have an integral representation
\begin{align}\begin{aligned}
 e^{-1}_\NS(r) + e^{-1}_\R(r) &= \zeta_0 \int\ud t e^{ \pi itr} \frac{e^{-\frac{1}{2}\pi tQ}}{G_\NS(Q+it)}  , \\
 e^{-1}_\NS(r) - e^{-1}_\R(r) &= \zeta_0 \int\ud t e^{ \pi itr} \frac{e^{-\frac{1}{2}\pi tQ}}{G_\R(Q+it)}  .
\end{aligned}\end{align}
The shift and reflection relations that they satisfy are as follows
\begin{align*}
e_\R\left(x-\frac{i \ub^{\pm1}}{2}\right) &= (1+ie^{\pi \ub^{\pm1}x})e_\NS\left(x+\frac{i\ub^{\pm1}}{2}\right), \\
e_\NS\left(x-\frac{i \ub^{\pm1}}{2}\right) &= (1-ie^{\pi \ub^{\pm1}x})e_\R\left(x+\frac{i\ub^{\pm1}}{2}\right), \\
e_\NS(x)e_\NS(-x) &= e^{-i\pi Q^2\slash8} e^{-i\pi(1-Q^2/2)\slash 3} e^{i\pi x^2\slash 2}, \\
e_\R(x)e_\R(-x) &= e^{i\pi\slash 2} e^{-i\pi Q^2\slash8} e^{-i\pi(1-Q^2/2)\slash 3} e^{i\pi x^2\slash 2}.
\end{align*}
Asymptotically, the functions $e_\NS$ and $e_\R$ behave as
\begin{align}
\label{asymptotic_superNS}
e_\NS(x) &= \left\{ \begin{array}{ll}
1 &\qquad  , x\to -\infty\\
e^{-i\pi Q^2\slash8} e^{-i\pi(1-Q^2/2)\slash 3} e^{i\pi x^2\slash 2} & \qquad , x\to +\infty\\
\end{array} \right.\\
e_\R(x) &= \left\{ \begin{array}{ll}
1 &, x\to -\infty\\
e^{i\pi\slash 2} e^{-i\pi Q^2\slash8} e^{-i\pi(1-Q^2/2)\slash 3} e^{i\pi x^2\slash 2} & ,x\to +\infty\\
\end{array} \right.
\end{align}
Also, we know that for self-adjoint operators $\mathsf{p},\mathsf{x}$ such that $[\mathsf{p},\mathsf{x}]=\frac{1}{ \pi i}$ they satisfy four pentagon relations
\begin{subequations}\label{superpentagonn_appendix}
	\begin{align}
	f_+(\mathsf{p}) f_+(\mathsf{x}) &=  f_+(\mathsf{x})f_+(\mathsf{x}+\mathsf{p})f_+(\mathsf{p}) -i f_-(\mathsf{x})f_-(\mathsf{x}+\mathsf{p})f_-(\mathsf{p}), \\
	f_+(\mathsf{p}) f_-(\mathsf{x}) &=  -i f_+(\mathsf{x})f_-(\mathsf{x}+\mathsf{p})f_-(\mathsf{p}) + f_-(\mathsf{x})f_+(\mathsf{x}+\mathsf{p})f_+(\mathsf{p}), \\
	f_-(\mathsf{p}) f_+(\mathsf{x}) &=  f_+(\mathsf{x}) f_+(\mathsf{x}+\mathsf{p}) f_-(\mathsf{p}) -i f_-(\mathsf{x})f_-(\mathsf{x}+\mathsf{p})f_+(\mathsf{p}) , \\
	f_-(\mathsf{p}) f_-(\mathsf{x}) &=  i f_+(\mathsf{x}) f_-(\mathsf{x}+\mathsf{p}) f_+(\mathsf{p}) - f_-(\mathsf{x})f_+(\mathsf{x}+\mathsf{p})f_-(\mathsf{p}),
	\end{align}
\end{subequations}
where $f_\pm(x) = e_\R(x) \pm e_\NS(x)$. The equations \eqref{superpentagonn_appendix} are equivalent to the supersymmetric analogues of the Ramanujan integral identities \cite{Hadasz:2007wi}
\begin{equation}\label{Ramanu}\begin{aligned}
\int_{-i\infty}^{i\infty} \frac{\ud\tau}{i} 
 e^{ \pi i \tau \beta} \left[ \frac{G_{\R}(\tau + \alpha)}{G_{\NS}
	(\tau + Q)} + \frac{G_{\NS}(\tau + \alpha)}{G_{\R}
	(\tau + Q)} \right] &= 2 \zeta_0^{-1} \frac{G_{\R}(\alpha)
	G_{\NS}(\beta)}{G_{\R}(\alpha+\beta)} , \\
\int_{-i\infty}^{i\infty} \frac{\ud\tau}{i} 
 e^{ \pi i \tau \beta} \left[ \frac{G_{\NS}(\tau + \alpha)}{G_{\NS}
	(\tau + Q)} + \frac{G_{\R}(\tau + \alpha)}{G_{\R}
	(\tau + Q)} \right] &= 2 \zeta_0^{-1} \frac{G_{\NS}(\alpha)
	G_{\NS}(\beta)}{G_{\NS}(\alpha+\beta)} , \\
\int_{-i\infty}^{i\infty} \frac{\ud\tau}{i} 
 e^{ \pi i \tau \beta} \left[  \frac{G_{\R}(\tau + \alpha)}{G_{\NS}
	(\tau + Q)} - \frac{G_{\NS}(\tau + \alpha)}{G_{\R}
	(\tau + Q)} \right] &= 2 \zeta_0^{-1} \frac{G_{\R}(\alpha)
	G_{\R}(\beta)}{G_{\NS}(\alpha+\beta)} , \\
\int_{-i\infty}^{i\infty} \frac{\ud\tau}{i} e^{ \pi i \tau \beta} \left[ \frac{G_{\NS}(\tau + \alpha)}{G_{\NS}
	(\tau + Q)} - \frac{G_{\R}(\tau + \alpha)}{G_{\R}
	(\tau + Q)} \right] &= 2 \zeta_0^{-1} \frac{G_{\NS}(\alpha)
	G_{\R}(\beta)}{G_{\NS}(\alpha+\beta)} .
\end{aligned}\end{equation}
Finally, we define
\begin{align}
 & g_\R(x) = e_\R(e^{\pi b x}), && g_\NS(x) = e_\NS(e^{\pi b x}) .
\end{align}

\subsection{Binomial and q-binomial identites}\label{Appendix.A.3}

The ordinary binomial and q-binomial formule have a continuous analogues. The continuous binomial formula for $x,y\in\mathbb{R}_{\geq0}$ is given by
\begin{equation}\begin{aligned}\label{continuousbinomialformulae}
(x+y)^{is} &= \int \frac{\ud t}{2\pi} {-is \choose -it}_\Gamma y^{it} x^{i(s-t)} , \\
(x-y)^{is} &= \int \frac{\ud t}{2\pi} {-is \choose -it}_\Gamma e^{\mp\pi t} y^{it} x^{i(s-t)} , \\
(-x-y)^{is} &= e^{\mp\pi s} \int \frac{\ud t}{2\pi} {-is \choose -it}_\Gamma y^{it} x^{i(s-t)} ,
\end{aligned}\end{equation}
where $s\in\mathbb{R}$, the sign depends on the choice of the branch of the logarithm and the continuous binomial coefficient is given explicitly in terms of gamma functions
\begin{equation}
 { s \choose t }_\Gamma = \frac{\Gamma(t)\Gamma(s-t)}{\Gamma(s)} .
\end{equation}
Moreover, for the non-commutative elements $u,v$, which satisfy $uv = q^2vu$ for $q=e^{i\pi \ub^2}$, one has the countinuous q-binomial formula
\begin{align}\label{continuous_q_binomialformula}
 &(u+v)^{it} = \ub \int \ud\tau {t \choose \tau}_\ub u^{i(t-\tau )} v^{i \tau },
\end{align}
where $t\in \mathbb{R}$ and a continuous version of the q-binomial coefficient has the form
\begin{equation}
 { t \choose \tau }_\ub = \frac{G_\ub(-\tau)G_\ub(-t+\tau)}{G_\ub(-t)} .
\end{equation}

\section{Toy model: continuous monomial algebra}
\label{appendixB}

As an instructive case on the road to the analysis of the continuous version of $U_q(sl(2))$ and $U_q(osp(1|2))$, one can consider a continuous version of the algebra of monomials, which behaves in a manner similar to the Cartan subalgebra of those Hopf algebras. The consideration of this algebra is much simpler, as it is both commutative and co-commutative. It provides an easy to understand example before the study of our main objects of interest in sections \ref{chapter3} and \ref{chapter4}. \\

In the case of an ordinary algebra of monomials, the algebra is spanned by a set of even elements $\{ y^n \}_{n\in\mathbb{N}}$ with the multiplication and co-multiplication relations as follows
\begin{align*}
& y^n y^m = y^{n+m}, \\
&\Delta(y^m) = (y_1+y_2)^n .
\end{align*}
Where it comes to the continuous case, let us start by considering a Hopf algebra $\mathcal{A}$ composed of the basis elements $\{ e(s,\epsilon) \}_{s\in\mathbb{R},\epsilon=\pm}$ given as follows
\begin{align}
e(s,+)&=\frac{1}{2\pi} \Gamma(-is)  e^{\frac{\pi s}{2}} (y)_{pv}^{is} ,\\
e(s,-)&=\frac{1}{2\pi}  \Gamma(-is)  e^{\frac{\pi s}{2}} (-y)_{pv}^{is} ,
\end{align}
where we use the following definitions of the principal value prescription for $y$
\begin{align*}
y_{pv}^{is} &={|y|}^{is}\Theta(y)+e^{-\pi s}|y|^{is}\Theta(-y), \\
(-y)_{pv}^{is} &=|y|^{is}\Theta(-y)+e^{-\pi s}|y|^{is}\Theta(+y),
\end{align*}
and where $\Theta$ is the step function. This definition is dictated by the fact that we do assume that $y$ is not a positive operator, and therefore its complex power needs to be made well-defined. 

For the continuous version of the algebra of monomials $\mathcal{A}$ the multiplication is as follows
\begin{align*}
e(s,+)e(s',+)&= \frac{1}{2\pi} { {-i(s+s')} \choose {-is} }_\Gamma e(s+s',+),\\
e(s,-)e(s',-)&= \frac{1}{2\pi} { {-i(s+s')} \choose {-is} }_\Gamma e(s+s',-) ,\\
e(s,+)e(s',-)&= \frac{1}{2\pi} { {-i(s+s')} \choose {-is} }_\Gamma \left[ \frac{e^{-\pi s'}(1-e^{-2\pi s})}{1-e^{-2\pi(s+s')}} e(s+s',+)+ \frac{e^{-\pi s}(1-e^{-2\pi s'})}{1-e^{-2\pi(s+s')}} e(s+s',-) \right] ,\\
&= e(s',-)e(s,+) ,
\end{align*}
and the co-multiplication
\begin{align*}
\Delta(e(s,+)) &= \int \ud \tau e(\tau,+) \otimes e(s-\tau,+),  \\
\Delta(e(s,-)) &= \int \ud \tau e(\tau,-) \otimes e(s-\tau,-) . 
\end{align*}

In addition to algebra $\mathcal{A}$ we consider also a dual Hopf algebra $\mathcal{A}^*$ composed of the basis elements $\{ \hat e(s,\epsilon) \}_{s\in\mathbb{R},\epsilon=\pm}$ defined by the equations
\begin{align}
\hat{e}(s,+) &= {|\hat y|}^{is} \Theta(\hat y) , \\
\hat{e}(s,-) &= {|\hat y|}^{is} \Theta(-\hat y) .
\end{align}
The multiplication of these elements are
\begin{align*}
\hat{e}(s,+)\hat{e}(s',+) &=\hat{e}(s+s',+), \\
\hat{e}(s,-)\hat{e}(s',-) &=\hat{e}(s+s',-) ,\\
\hat{e}(s,-)\hat{e}(s',+) &= 0 , \\
\hat{e}(s,+)\hat{e}(s',-) &= 0 ,
\end{align*}
and co-multiplication
\begin{align*}
{\hat \Delta}(\hat e(s,+))&=\int \frac{\ud\tau}{2\pi} {-is \choose -it}_\Gamma \left[ \hat e(\tau,+)\otimes\hat e(s-\tau,+) + \hat e(\tau,+)\otimes\hat e(s-\tau,-) \frac{e^{-\pi (s-\tau)} (1-e^{-2\pi\tau})}{1-e^{-2\pi s}} + \right. \\
&\left. +
\hat e(\tau,-)\otimes\hat e(s-\tau,+) \frac{e^{-\pi\tau}(1-e^{-2\pi(s-\tau)})}{1-e^{-2\pi s}} \right],\\
{\hat \Delta}(\hat e(s,-))&=\int \frac{\ud\tau}{2\pi} {-is \choose -it}_\Gamma \left[ \hat e(\tau,-)\otimes\hat e(s-\tau,-) +  \hat e(\tau,+)\otimes\hat e(s-\tau,-) \frac{e^{-\pi\tau}(1-e^{-2\pi(s-\tau)})}{1-e^{-2\pi s}} + \right. \\
&\left. +
\hat e(\tau,-)\otimes\hat e(s-\tau,+) \frac{e^{-\pi (s-\tau)} (1-e^{-2\pi\tau})}{1-e^{-2\pi s}} \right] .
\end{align*}
By observation one can clearly see that the coefficients of multiplication and co-multiplication for the algebras $\mathcal{A}$ and $\mathcal{A}^*$ defined by the following equations
\begin{align*}
e(s,\epsilon) e(s',\epsilon') &= \sum_{\epsilon''} \int_{\mathbb{R}} \ud\sigma \, m_{s,\epsilon;s',\epsilon'}^{\sigma,\epsilon''} e(\sigma,\epsilon''), \\
\Delta(e(s,\epsilon)) &= \sum_{\epsilon',\epsilon''} \int_{\mathbb{R}^2} \ud\sigma\ud\sigma' \, \mu_{s,\epsilon}^{\sigma,\epsilon';\sigma',\epsilon''} e(\sigma,\epsilon')\otimes e(\sigma',\epsilon''), \\
\hat e(s,\epsilon) \hat e(s',\epsilon') &= \sum_{\epsilon''} \int_{\mathbb{R}} \ud\sigma \, {\hat m}^{s,\epsilon;s',\epsilon'}_{\sigma,\epsilon''} \hat e(\sigma,\epsilon''), \\
{\hat \Delta}(\hat e(s,\epsilon)) &= \sum_{\epsilon',\epsilon''} \int_{\mathbb{R}^2} \ud\sigma\ud\sigma' \, {\hat \mu}^{s,\epsilon}_{\sigma,\epsilon';\sigma',\epsilon''} \hat e(\sigma,\epsilon')\otimes \hat e(\sigma',\epsilon'') ,
\end{align*}
satisfy the relations
\begin{align*}
m_{s,\epsilon_1;s',\epsilon_2}^{\sigma,\epsilon_3} &= \hat{\mu}_{s,\epsilon_1;s',\epsilon_2}^{\sigma,\epsilon_3}, \\
\hat{m}^{s,\epsilon_1;s',\epsilon_2}_{\sigma,\epsilon_3} &= \mu^{s,\epsilon_1;s',\epsilon_2}_{\sigma,\epsilon_3} ,
\end{align*}
i.e. we see that those two algebras are indeed dual to each other. This allows us to define the multiplication relations between the elements of the Heisenberg algebra $H(\mathcal{A})$
\begin{align*}
(1 \otimes e(s,\pm))( \hat e(s',\pm) \otimes 1) &= \int \frac{\ud\tau}{2\pi} { {-i s'} \choose {-i\tau}}_\Gamma \hat e(s'-\tau,\pm) \otimes e(s-\tau,\pm), \\
(1 \otimes e(s,\pm)) ( \hat e(s',\mp) \otimes 1) &= 0 .
\end{align*}
If we analytically continue those expressions to $s,s'\in -i\mathbb{N}_{\geq0}$ we can find that the generators satisfy the canonical commutation relations
\begin{equation}
 y \hat y = \hat y y + \frac{1}{i} ,
\end{equation}
where we denoted $1\otimes y$ as $y$ and ${\hat y} \otimes 1$ as ${\hat y}$ for the sake of brevity. Using the realisation in terms of generators, we can compute the canonical element $S$ using the formula \eqref{canonicalelementdef}, which gives
\begin{equation}
S = e^{i y\otimes\hat y} ,
\end{equation}
to compute which we used the Mellin transformation for an exponential 
\begin{align}\label{mellin_transform_exponential}
e^{\pm ix} &= \int_{-\infty}^\infty \frac{\ud t}{2\pi} \Gamma(-it) e^{\mp\frac{\pi t}{2}} x^{it} ,
\end{align}
where $x>0$.

\end{document}